%% file: tziolas.tex
\documentclass{amsart}
\usepackage{amsthm}
\usepackage{amsmath}
\usepackage{amstext}
\usepackage{amsfonts}
\usepackage{latexsym}
\usepackage{amscd}
\usepackage{xypic}
\newcommand{\la}{\ensuremath{\longrightarrow}}

\newcommand{\pone}{\ensuremath{\mathbb{P}^{1}}}

\newcommand{\sheaf}{\ensuremath{\mathcal{O}}}

\newcommand{\ext}{\ensuremath{\underline{Ext}}}

\newtheorem{theorem}{Theorem}[section]
\newtheorem{lemma}[theorem]{Lemma}
{\theoremstyle{definition} \newtheorem{definition}[theorem]{Definition}}
{\theoremstyle{definition} }
{\theoremstyle{definition} \newtheorem{example}[theorem]{Example}}
{\theoremstyle{definition} \newtheorem{remark}[theorem]{Remark}}
\newtheorem{proposition}[theorem]{Proposition}
\newtheorem{corollary}[theorem]{Corollary}

\begin{document}

\title[$\mathbb{Q}$-Gorenstein deformations of nonnormal surfaces]{$\mathbb{Q}$-Gorenstein deformations of nonnormal surfaces}
\author{Nikolaos Tziolas}
\address{Department of Mathematics, University of Cyprus, P.O. Box 20537, Nicosia, 1678, Cyprus}
\email{tziolas@ucy.ac.cy}
\thanks{The author was partially supported by a JSPS fellowship, \# PE 04054}

%    General info
\subjclass{Primary 14E30, 14E35}
%\date{April 21, 1999}

%\dedicatory{This paper is dedicated to our authors.}

\keywords{Algebraic geometry}

\begin{abstract}
Let $\Delta \subset H$ be the germ of a non-normal surface along a proper curve with smooth components such that the high index points of $H$ are semi-log-terminal 
and the Gorenstein singular points are semi-log-canonical. We describe the sheaf $T^1_{qG}(H)$ of $\mathbb{Q}$-Gorenstein deformations of $H$ and we show that 
the divisorial part of the support of $T^1_{qG}(H)$ is $\Delta$. Moreover, we show that if $C$ is any irreducible component of $\Delta$, then 
the locally free part of $T^1_{qG}(H)\otimes \sheaf_C$ is a line bundle on $C$ 
and we obtain a formula for its degree. Finally we obtain criteria 
for $\Delta\subset H$ to have $\mathbb{Q}$-Gorenstein terminal smoothings.   
\end{abstract}

\maketitle

\input{introduction-tziolas}
\input{section1-tziolas}

\input{section2-tziolas}
\input{section3-tziolas}

\input{bibliography-tziolas}
\end{document}

%% file: introduction-tziolas.tex
\section{Introduction}

The purpose of this paper is to describe the $\mathbb{Q}$-Gorenstein smoothings of a germ $\Delta \subset H$ of a nonnormal surface with semi-log-canonical (\textit{slc}) singularities along a 
proper curve with smooth irreducible components. $\mathbb{Q}$-Gorenstein smoothings of such surface germs provide important information about the components of the versal deformation space 
$\mathrm{Def}(X)$ of an isolated rational surface singularity $0\in X$. They are also closely related to 3-fold terminal extremal neighborhoods, a fundamental class of birational maps that 
appear in the three dimensional minimal model program, and they help to understand the boundary of the compactification of the moduli space of surfaces of general type. 

A stable surface is a proper two-dimensional reduced scheme $H$ such that $H$ has only semi-log-canonical singularities and $\omega_H^{[k]}$ is locally free and ample for some $k>0$. 
The boundary of the compactification of surfaces of general type consists of smoothable stable surfaces. Therefore it is interesting to know which stable surfaces are smoothable.

Let $0\in X$ be the germ of a rational surface singularity. Koll\'{a}r~\cite{Ko91} has made a series of conjectures concerning the components of the versal deformation space $\mathrm{Def}(X)$. 
A $P$-resolution of of $X$~\cite{KoBa88} is a proper birational map $g \colon H\la X$ such that 
\begin{enumerate}
\item $H$ is Cohen Macaulay and Gorenstein outside finitely many points
\item $K_H$ is $g$-ample
\item $H$ has a $\mathbb{Q}$-Gorenstein smoothing. 
\end{enumerate}
Koll\'{a}r conjectures that there is a one-to-one correspondence between the components of the versal deformation space $\mathrm{Def}(X)$ of $X$, and $qG$-components of deformations 
of $P$-resolutions~\cite{Ko91}. This conjecture is known to be true for Gorenstein surface singularities, quotient singularities~\cite{KoBa88} and rational quadruple points~\cite{Ste91}. 

If a $P$-resolution $H$ of a rational surface singularity is normal, then its deformation theory can be described by the local deformations around the singular points. However, 
there are examples~\cite{Ko91} that show that a $P$-resolution may not be normal. In such a case it is possible that locally $H$ is smoothable but $H$ itself is not. Therefore 
it is important to give a more detailed study of the deformation theory of nonnormal surface germs along smooth proper curves.

There is also a close relation between $\mathbb{Q}$-Gorenstein smoothings of surface germs $\Delta \subset H$ as above and 3-fold terminal extremal neighborhoods. 
A 3-fold terminal extremal neighborhood~\cite{Ko-Mo92} is a proper birational map $\Delta \subset Y \stackrel{f}{\la} X \ni P$ such that $Y$ is the germ of a 3-fold along a proper curve $\Delta$, 
$\Delta_{red}=f^{-1}(P)$, $X$ and $Y$ are terminal, and $-K_Y$ is $f$-ample. Then $Y$ is a one-parameter $\mathbb{Q}$-Gorenstein smoothing of the general member $H \in |\sheaf_Y|$. In 
general $H$ may or may not be normal and its singularities may be hard to describe. Koll\`{a}r and Mori~\cite{Ko-Mo92} showed that given a terminal extremal neighborhood as above, the general 
members of $|-K_Y|$ and $|-K_X|$ have DuVal singularities. If the neighborhood is isolated and the general member of $|-K_X|$ is of type $A_n$, then the general member $H \in |\sheaf_Y|$ 
has semi-log-canonical (\textit{slc}) singularities, and in fact the points of index bigger than 1 are semi-log-terminal (\textit{slt})~\cite{Tzi04}. The same is true 
if the neighborhood is divisorial and the general member $S\in |-K_X|$ that contains $\Gamma=f(E)$, where $E$ is the $f$-exceptional divisor, is of type $A_m$. Hence a detailed description 
of the $\mathbb{Q}$-Gorenstein deformations of non-normal surface germs $\Delta \subset H$ with \textit{slc} singularities will provide information about the structure of terminal extremal neighborhoods. 

Let $\Delta\subset H$ be the germ of a surface with \textit{slc} singularities along a smooth proper curve such that the 
points of index bigger than 1 are \text{slt}. By the classification of \textit{slc} singularities~\cite{KoBa88}, if such a germ has a terminal smoothing then its singularities are either 
normal crossing points, pinch points, degenerate cusps of embedding dimension at most 4 or \textit{slt} points analytically isomorphic to $(xy=0)/\mathbb{Z}_n(a,-a,1)$, $(a,n)=1$. 
Such singularities we call singularities of class \textit{qG}. A surface germ $\Delta \subset H$ with singularities of class \textit{qG} has local $\mathbb{Q}$-Gorenstein smoothings 
but there is not necessarily a global one extending the local. 

The purpose of this paper is to study when the local $\mathbb{Q}$-Gorenstein smoothings extend to a global. 
To put the problem in the proper setting, 
we define deformation functors $Def^{qG}(H)$ and $Def^{qG}_{loc}(H)$ that parametrize global and collections of local $\mathbb{Q}$-Gorenstein deformations, and we study the natural 
transformation $\phi \colon Def^{qG}(H) \la Def^{qG}_{loc}(H)$ and in particular the induced map on tangent spaces. 

Theorem~\ref{main-theorem} describes the tangent space of $Def^{qG}_{loc}(H)$. This is $H^0(T^1_{qG}(H))$, 
where $T^1_{qG}(H)\subset T^1(H)$ is the sheaf of local $\mathbb{Q}$-Gorenstein deformations of $H$. It is shown that the divisorial part of 
the support of $T^1_{qG}(H)$ is $\Delta$. Let $C$ be an irreducible component of $\Delta$ and $L_C$ the locally free part of $T^1_{qG}(H)\otimes\sheaf_C$. It is 
shown that this is a line bundle and a formula for its degree 
is given which involves the self intersection $\tilde{C}^2$ and some analytic 
invariants of the singularities of $H$, where $\tilde{C}=\pi^{-1}(C)$ and $\pi \colon \tilde{H} \la H$ is the normalization of $H$. At this point I must mention 
that a special case of Theorem~\ref{main-theorem} appeared without proof in~\cite{Ko91}, and in fact this was the original motivation for this work.  

If $H^2(T_H)=0$, then Proposition~\ref{def-map} says that the natural transformation of the deformation functors $Def^{qG}(H) \la Def^{qG}_{loc}(H)$ is smooth. In particular this is the case 
when $H$ is a modification of an isolated singularity. Moreover, if every irreducible component of $\Delta$ is a smooth rational curve, then $H$ has global $\mathbb{Q}$-Gorenstein smoothings 
if and only if $d_c = \deg L_C \geq 0$.

As applications, Corollary~\ref{global-def} gives a necessary and sufficient condition for a germ $\Delta \subset H$ 
to have a one parameter terminal $\mathbb{Q}$-Gorenstein smoothing in the case that $\Delta$ is a rational cycle of curves, and Corollary~\ref{extremal} shows that under some conditions there exists a 
terminal extremal neighborhood $\Delta \subset Y \stackrel{f}{\la} X \ni P$, such that $H \in |\sheaf_Y|$.

I would like to thank J\'anos Koll\'ar for many fruitfull discussions and exchange of ideas. Also Lemma~\ref{patching} was communicated to  me by him.

%% file: section1-tziolas.tex
\section{Preliminaries}
Let $\mathcal{F}$ be a coherent sheaf on a scheme $X$. We denote $\mathcal{F}^{[n]}=({\mathcal{F}^{\otimes n}})^{\ast\ast}$. 
\begin{definition}
Let $X$ be either a variety, an analytic space or the germ of a singularity, such that it is Cohen Macaulay and Gorenstein in codimension 1.
Let $Y$ be the total space of a one parameter deformation of $X$. Then we say that $Y$ is a $\mathbb{Q}$-Gorenstein deformation if $\omega_Y^{[n]}$ is locally 
free for some $n$. The assumptions on $X$ assure that this makes sense even for nonnormal varieties. 
\end{definition}
We will now define three deformation functors. Let $Art(\mathbb{C})$ be the category of finite local Artin $\mathbb{C}$-algebras and $Sh(X)$ the category of sheaves of sets on $X$.
\begin{definition}[Definition 3.17~\cite{KoBa88}]          
The functor of global $\mathbb{Q}$-Gorenstein deformations is the functor $Def^{qG}(X) \colon Art(\mathbb{C}) \la Sets$ such that for any finite local $\mathbb{C}$-algebra $A$, 
 $Def^{qG}(X)(A)$ is the set of isomorphism classes of flat morphisms $f\colon Y \la S =\mathrm{Spec}(A)$, such that $Y\otimes k(A) \cong X$, the sheaf 
$\omega^{[n]}_{Y/S}$ is invertible for some $n$.
\end{definition}
It is not immediately clear that $Def^{qG}(X)$ as defined above is indeed a functor. For this to be true, the property that the relative dualizing sheaf $\omega_{Y/S}$ 
is $\mathbb{Q}$-Gorenstein should be stable under base
extension. However, this is true~\cite[Lemma 2.6]{Has-Kov04} and therefore  $Def^{qG}(X)$ is a functor 
and moreover a subfunctor of the versal deformation functor $Def(X)$.

There are three cases when it is known that $Def(X)$ is pro-represented by an analytic space $\mathrm{Def}(X)$.
\begin{enumerate}
\item If $(P\in X)$ is the germ of an isolated singularity~\cite{Gra72}~\cite{Ko-Mo92},.
\item If $X$ is a compact complex space~\cite{Gra74}~\cite{Ko-Mo92}
\item If there is a proper birational morphism $f \colon X \la Y$, where $Y$ is the germ of an isolated singularity and $R^1f_{\ast}\sheaf_X =0$~\cite{Wa76}~\cite{Ko-Mo92}.
\end{enumerate}
In general it is not known whether $Def^{qG}(X)$ is pro-representable in the above cases. However, if $X$ is normal, $\dim X =2$ and its singularities of index bigger than 1 are all of 
class $T$~\cite{KoBa88}, then $Def^{qG}(X)$ is pro-represented by a closed subspace $\mathrm{Def}^{qG}(X)$ of $\mathrm{Def}(X)$. 
In fact,  $\mathrm{Def}^{qG}(X)= \pi^{-1} (\prod _P \mathrm{Def}^{qG} (P\in X))$, where $\pi$ is the natural map 
\[
\mathrm{Def}(X) \la \prod _P \mathrm{Def} (P\in X)
\]
$P \in X$ is a singular point of index bigger than 1 and $\mathrm{Def}^{qG} (P\in X)$ is the $qG$-Gorenstein component of $\mathrm{Def}(P\in X)$ 
which is known to exist~\cite{KoBa88}.

The purpose of this paper is to describe the tangent space of the functor $Def^{qG}(X)$ when $X$ has singularities of class $qG$~[Definition~\ref{qG}] 
and to obtain criteria for the existence of 
$\mathbb{Q}$-Gorenstein smoothings of $X$. The problem of the pro-representability of $Def^{qG}(X)$ will not be touched here. 
\begin{definition}
The functor of $\underline{Def}^{qG}(X)$ is the functor \[
\underline{Def}^{qG}(X) \colon Art(\mathbb{C}) \la Sh(X)\]
 defined as follows. For any finite local 
$\mathbb{C}$-algebra $A$,  $\underline{Def}^{qG}(X)(A)$ is the sheaf defined by $\underline{Def}^{qG}(X)(A)(U)=Def^{qG}(U)(A)$, for any affine open subset $U$ of $X$.
\end{definition}
Note that $\underline{Def}^{qG}(X)(A)$ is the sheaf associated to the presheaf $F$ defined by $F(V) = Def^{qG}(V)(A)$ for any open set $V$.
 
Let $k$ be a field and $F\colon Art(k) \la \mathcal{C}$ be any functor from the category of finite local $k$-algebras to a category $\mathcal{C}$. Then $F(k[t]/(t^2))$ is called 
the tangent space of $F$~\cite{Sch68}.
\begin{definition}
We denote by $ \mathbb{T}^1_{qG}(X)$, $T^1_{qG}(X)$ the tangent spaces of the functors $Def^{qG}(X)$ and  $\underline{Def}^{qG}(X)$, respectively.  $ \mathbb{T}^1_{qG}(X)$ is a finite dimensional 
vector space over $\mathbb{C}$ and  $T^1_{qG}(X)$ is a sheaf of $\sheaf_X$-modules defined as follows. 
For any affine open subset $U \subset X$, $T^1_{qG}(X)(U)$ is the $\sheaf_X(U)$-module
of isomorphism classes of first order $\mathbb{Q}$-Gorenstein deformations of $U$. If $U$ is Gorenstein, then $T^1_{qG}(X)(U)=T^1(U)$, the space of first order deformations of $U$.  
\end{definition} 

Standard results from deformation theory show that there is an exact sequence
\[
0 \la H^1(T_X) \la  \mathbb{T}^1_{qG}(X) \la H^0(T^1_{qG}(X)) \la H^2(T_X)
\]
and therefore in order to describe the first order $\mathbb{Q}$-Gorenstein deformations of $X$, and in particular the existence of smoothings,
one needs to study the sheaf $T^1_{qG}(X)$ and the obstruction space $H^2(T_X)$.The purpose of this paper is to describe $T^1_{qG}(X)$.
\begin{definition}
The functor of local $\mathbb{Q}$-Gorenstein deformations is the functor $
Def^{qG}_{loc}(X) \colon Art(\mathbb{C}) \la Sets$ defined by \[
Def_{loc}^{qG}(X)(A)=H^0(\underline{Def}^{qG}(X)(A))\]
\end{definition}
From the definition of the above functors it is clear that there exists a natural transformation $Def^{qG}(X) \la Def^{qG}_{loc}(X)$. If $X$ satisfies any of the 
conditions (1), (2) or (3), above, then both functors are pro-represented by $\mathrm{Def}^{qG}(X)$ and $\mathrm{Def}^{qG}_{loc}(X)$. 
$\mathrm{Def}^{qG}_{loc}(X)$ parametrizes collections of local deformations and its tangent space is $H^0(T^1_{qG}(X))$. 

In this paper we are interested to know when local deformations exist globally. In order to do this we will study $Def^{qG}_{loc}(X)$, and in particular its tangent space $H^0(T^1_{qG}(X))$, 
and we will find cases when the natural morphism 
$\phi \colon Def^{qG}(X) \la Def^{qG}_{loc}(X)$ is smooth.

\begin{definition}
A surface singularity $(P\in X)$ is called a normal crossing point (resp. pinch point) if it is analytically isomorphic to $(xy=0)\subset \mathbb{C}^3$ (resp. $(x^2-y^2z=0)\subset \mathbb{C}^3$).
\end{definition}
\begin{definition}
Let $(P \in X)$ be the germ of a nonnormal surface singularity. A map $f \colon Y \la X$ is called a semi-resolution of $X$ if the following conditions are satisfied:
\begin{enumerate}
\item $f$ is proper
\item $Y$ is semismooth
\item If $D_Y$ is the double curve of $Y$, then no component of $D_Y$ is mapped to a point by $f$
\item There is a finite set $S \subset X$ such that $f \colon f^{-1}(X-S) \la X-S$ is an isomorphism.
\end{enumerate}
\end{definition}
\begin{definition}
Let $P\in X$ be a $\mathbb{Q}$-Gorenstein surface singularity such that $X-P$ is semismooth. Let $f\colon Y \la X$ be a good semi-resolution of $X$. We can write 
$\omega^n_Y=f^{\ast}\omega^{[n]}_X\otimes \sheaf_Y (\sum a_iE_i)$, where $E_i$ are $f$-exceptional. Then $(P\in X)$ is called
\begin{enumerate}
\item semi-log-terminal(\textit{slt}) if $a_i >-1$, for all $i$,
\item semi-log-canonical(\textit{slc}) if $a_i\geq -1$, for all $i$,
\end{enumerate}
\end{definition}
\begin{definition}[~\cite{Ba83}]
Let $X$ be a Gorenstein surface singularity which has a minimal semi-resolution $f \colon Y \la X$. $X$ is called a degenerate cusp if $X$ is not normal and the exceptional divisor 
is a cycle of smooth rational curves or a rational nodal curve.
\end{definition}
It is known~\cite{KoBa88} that a Gorenstein surface singularity $P\in X$ is semi-log-canonical, iff $X$ is either simple elliptic, a cusp, a degenerate cusp or semi-canonical. 
Moreover, from the classification of \textit{slt} singularities~\cite{KoBa88} it follows that the only \textit{slt} singularities that admit $\mathbb{Q}$-Gorenstein smoothings are analytically 
isomorphic to $(xy=0)/\mathbb{Z}_n(a,-a,1)$, with $(a,n)=1$, and their $\mathbb{Q}$-Gorenstein deformations are given by $(xy+tf(t,z^n)=0)/\mathbb{Z}_n(a,-a,1,0)$. 
If a semi-log-canonical surface germ $X$ has a terminal $\mathbb{Q}$-Gorenstein smoothing, then the degenerate cusps of $X$ must have embedding dimension at most 4. This motivates 
the following.

\begin{definition}\label{qG}
A nonnormal surface singularity $P\in X$ is called a singularity of class $qG$, if it is analytically isomorphic to one of the following
\begin{enumerate}
\item Normal crossing point: $(xy=0)\subset \mathbb{C}^3$
\item Pinch point: $(x^2-y^2z=0)\subset \mathbb{C}^3$
\item Degenerate cusp of embedding dimension at most 4
\item Semi-log-terminal (\textit{slt}) : $(xy=0)/\mathbb{Z}_n(a,-a,1)$, $(a,n)=1$.
\end{enumerate} 
\end{definition}  

Let $(P\in X)$ be a degenerate cusp and $f \colon Y\la X$ the minimal semi-resolution. Let $\Gamma$ be the reduced exceptional locus. The next lemma is a classification of degenerate cusps 
of embedding dimension at most 4.

\begin{lemma}[Karras~\cite{Ka77}, Shepherd-Barron~\cite{Ba83}, Stevens~\cite{Ste98}]\label{equations}
Let $P\in X$ be a degenerate cusp and let $f\colon Y \la X$ be the minimal semi-resolution. Let $\Gamma = f^{-1}(P)_{{\mathrm{red}}}$. Then 
\begin{enumerate}
\item \[
\mathrm{mult}_P(X) =\mathrm{max}\{2,-\Gamma^2\} \]
\item \[
\mathrm{embdim}_P(X)=\mathrm{max}\{3,-\Gamma^2\} \]
\item If $\Gamma^2=-1$, then in local analytic coordinates \[
(P\in X)\cong (x^2=y^3+y^2z^2)\subset \mathbb{C}^3 \]
\item If $\Gamma^2=-2$ then in suitable local analytic coordinates
\[
(P\in X)\cong (x^2+z^2(z^{n+1}-y^2)=0)\subset \mathbb{C}^3 \]
with $2 \leq n \leq \infty$. In the case that $n=\infty$ we set $z^n =0$. 
If $n<\infty$, then the singular locus of $X$ is a smooth irreducible curve. If $n=\infty$, then the singular locus of $X$ is a reducible curve 
with exactly two smooth irreducible components. We will call such a singularity a degenerate cusp of type $T^2_n$.
\item If $\Gamma^2 = -3$ then in local analytic coordinates \[
(P \in X) \cong (x^{p+2} +y^{q+2} -xyz=0) \subset \mathbb{C}^3 \]
where $1\leq p,\; q \leq \infty$ and if $u$ is any of $x,\; y,\; z$, we set $u^{\infty}=0$. Moreover, the singular locus of $X$ has exactly $1+k$ smooth irreducible components, 
where $k$ is the number of exponents that are $\infty$. We will call such a singularity  a degenerate cusp of type $T^3_{p,q}$.
\item If $\Gamma^2=-4$ then in local analytic coordinates \[
(P\in X)\cong (xy-z^p-t^q=0,\; zt-x^r=0)\subset \mathbb{C}^4\]
where $2 \leq p,\; q,\; r \leq \infty$ and if $u$ is any of $x,\; y,\; z,\; t$, we set $u^{\infty}=0$. Moreover, the singular locus of $X$ has exactly $1+k$ smooth irreducible components, 
where $k$ is the number of exponents that are $\infty$.
We will call such a singularity a degenerate cusp of type $T^4_{p,q,r}$.
\end{enumerate}
\end{lemma}
In particular we see that any degenerate cusp of embedding dimension at most 4 is a complete intersection and in fact a degenerate cusp is a complete intersection if and only if it 
has embedding dimension at most 4~\cite{Ste98}. Hence they all have smoothings and in fact $T^1_{qG}(P\in X)=T^1(P\in X)$. 

Let $P\in H$ be a degenerate cusp of embedding dimension $4$ with reducible singular locus $\Delta$( this is the case when at least one of $p$, $q$, $r$ is infinity), 
and let $C$ be a component of $\Delta$. There are many isomorphic ways that one can choose local coordinates to describe the equations $C \subset H$. For example, 
let $P\in H$ be of type $T^4_{p,\infty,\infty}$ and $C$ the curve $x=y=z=0$. This is isomorphic to the curve $z=t=x=0$ in a degenerate cusp of type $T^4_{\infty,\infty,p}$. 
In order to get a description of $T^1_{qG}(H)$ we will have to make a choice of local coordinates. This is mostly due to the reason that we want to describe the invariants that 
appear in~\ref{main-theorem} and correspond to degenerate cusps with reducible singular locus, as limiting cases of the ones that appear in the case of degenerate cusps with irreducible 
singular locus. The choices that we make in the following corollary are purely for technical reasons.

From Lemma~\ref{equations} it follows that 
\begin{corollary} 
Let $P\in H$ be a 
degenerate cusp of embedding dimension $4$ and with reducible singular locus $\Delta$. Let $C$ be a component of $\Delta$. Then in suitable local analytic coordinates 
\begin{enumerate}
\item $C=(x=z=t=0)$ and $H=(xy-z^p-t^q=zt=0)\subset \mathbb{C}^4$, with $p,\; q <\infty$. In this case we will call the germ $C \subset H$ of type $T^4_{p,q,\infty}$.
\item $C=(x=z=t=0)$ and $H=(xy-z^p=zt-x^r=0)\subset \mathbb{C}^4$, with $p,\; r <\infty$. In this case we will call the germ $C\subset H$ of type $T^4_{p,\infty,r}$.
\item $C=(x=y=z=0)$ and $H=(xy-z^p=zt=0)\subset \mathbb{C}^4$, with $p < \infty$. In this case we will call the germ $C\subset H$ of type $T^4_{p,\infty,\infty}$.
\item  $C=(x=y=z=0)$ and $H=(xy=zt-x^r=0)\subset \mathbb{C}^4$, with $r < \infty$. In this case we will call the germ $C\subset H$ of type $T^4_{\infty,\infty,p}$.
\item $C=(x=y=z=0)$ and $H=(xy=zt=0)\subset \mathbb{C}^4$. In this case we will call the germ $C\subset H$ of type $T^4_{\infty,\infty,\infty}$.
\end{enumerate}
\end{corollary}

\begin{definition}
Let $(P\in X)$ be a surface cyclic quotient singularity of type $1/n(1,a)$. Let $U $ be the minimal resolution of $H$, $E_1,\ldots,E_m$ the exceptional divisors and 
\[
\stackrel{E_1}{\circ}\mbox{\noindent ---} \stackrel{E_2}{\circ}\mbox{\noindent ---}
\cdots \mbox{\noindent ---}\stackrel{E_m}{\circ} \]
be the fundamental cycle. Let $C$ be a smooth curve in $X$ going through the singularity. Then $C$ will be called of type $FC_k$ if the birational transform $C^{\prime}$ of $C$ 
in $U$ intersects $E_k$.
\end{definition}

%% file: section2-tziolas.tex
\section{Local to Global}
The next Theorem describes the tangent space of $\underline{Def}^{qG}(H)$ and $Def^{qG}_{loc}(H)$.
\begin{theorem}\label{main-theorem}
Let $\Delta \subset H$ be the germ of a surface along a reduced proper curve. Assume that every irreducible component of $\Delta$ is smooth and that locally
along $\Delta$ the surface $H$ has singularities of class $qG$.
Let $\pi\colon \tilde{H} \la H$ be the normalization of $H$. Let $C$ be an irreducible component of $\Delta$ and $\tilde{C}$ the divisorial part of $\pi^{-1}(C)$.
Let $p$ be the number of pinch points and $c_i$ the number of degenerate cusps with $\Gamma^2=-i$, $1\leq i \leq 4$ of $H$ on $C$.
Let $U_3$, $U_4$ be the sets of points $P$ of $H$ that  lie on $C$ such that $P\in H$ is a degenerate cusp with $\Gamma^2=-3$ and $ -4$ respectively. Then
\begin{enumerate}
\item The support of $T^1_{qG}(H)$ is 1-dimensional. Its divisorial part is $\Delta$ and it has an embedded point over any pinch point or degenerate cusp
of type other than  $T^3_{\infty, \infty}$ and $T^4_{\infty,\infty,\infty}$. Moreover, let $C$ be any component of $\Delta$. Then,
\[
T^1_{qG}(H) \otimes \sheaf_C = L_C \oplus (\oplus_{i}\mathbb{C}_{{P_i}}) \]
where $L_C$ is a line bundle on $C$ and $\mathbb{C}_{{P_i}}$ is a 1-dimensional torsion sheaf supported at a degenerate cusp of embedding dimension 4 that lies on $C$.
\item \[
\deg L_C ={\tilde{C}}^2 +p+2c_1+2c_2+\sum_{P\in U_3} \alpha_3(P) +\sum_{Q\in U_4} \alpha_4(Q)
\]
where \[
\alpha_3(P)=1+\frac{p+q}{pq+p+q} \]
if $P\in U_3$ and $P\in C \subset H$ is of type $T^3_{p,q}$, and \[
\alpha_4(Q)=\frac{r(p+q)-4}{rpq-p-q}\]
if $Q\in U_4$ and $P\in C \subset H$ is of type $T^4_{p,q,r}$. If any of $p$, $q$, $r$ is infinity, then we understand $\alpha_3(P)$ and $\alpha_4(Q)$ to mean the corresponding limits at infinity.
\end{enumerate}
\end{theorem}

The proof of Theorem~\ref{main-theorem} is given in section 4.

\noindent\textbf{Remark:} From the equations defining a degenerate cusp, we see that a singularity $P\in H$ of type $T^3_{p,q}$ can be considered as of type $T^4_{p+2,q+2,1}$.
Then the contribution to the degree of $L_C$ in 3.2 above is $\alpha_3(P)$. One may expect that this should be the same as $\alpha_4(P)$. However, $\alpha_3(P)=\alpha_4(P)+1$.
The reason that there is an extra 1 is that we are assuming that $r\geq 2$ and in this case 3.1 shows that
$L_C$ has a 1-dimensional torsion at $P$ that accounts for the extra 1.

\begin{proposition}\label{def-map}
Let $\Delta \subset H$ be a surface germ as in Theorem~\ref{main-theorem}, and $C_i$, $i=1, \ldots ,k$ the irreducible components of $\Delta$. Let 
$L_{C_i} \in \mathrm{Pic}(C_i)$ the locally free part of $T^1_{qG}(H) \otimes \sheaf_{C_i}$. Assume that $H^2(T_H)=0$. Then
\begin{enumerate}
\item The natural map of functors \[
\phi \colon Def^{qG}(H) \la Def^{qG}_{loc}(H) \]
is smooth.
\item Assume that $C_i \cong \mathbb{P}^1$, $\forall i$. Then a one-parameter $\mathbb{Q}$-Gorenstein terminal smoothing $X$ of $H$ exists if and only if $d_i= \deg( L_{{C_i}})\geq 0$. $\forall i$. 
Moreover, in this case $X$ can be chosen so that it has isolated cyclic quotient singularities of type $1/n(a,-a,1)$ and ordinary double points $xy-zt=0$.
\end{enumerate}
\end{proposition}
\noindent\textbf{Remark:} In general it is rather difficult to calculate $H^2(T_H)$. 
However if $H$ is a modification of an isolated surface singularity $0\in S$, i.e., there is a birational morphism 
$f\colon \Delta \subset H \la S \ni 0$, then by the formal functions theorem it follows that $H^2(T_H)=0$. 
Moreover, if $R^1f_{\ast}\sheaf_H =0$ (for example if the singularity $0\in S$ is rational), then $\Delta$ has rational components and hence 
the assumptions of Proposition~\ref{def-map} are satisfied.

\begin{proof}
The map $\phi$ is a smooth map of functors iff \[
Def^{qG}(H)(B) \la Def^{qG}(H)(A)\times_{{Def^{qG}_{loc}(H)(A)}} Def^{qG}_{loc}(H)(B) \]
is surjective for any small extension of finite local $\mathbb{C}$-algebras~\cite{Li-Sch67}\[
0 \la N \la B \la A \la 0\]
In fact it is sufficient to consider only the case when $N=A/m_A=\mathbb{C}$.

In general, let $X$ be a scheme and $\mathcal{F}$ a coherent sheaf on it. Then the space $Ex(X,\mathcal{F})$ of infinitesimal extensions of $X$ by $\mathcal{F}$ fits in the
following exact sequence
\begin{equation}
0 \la H^1(\underline{Hom}(\Omega_X,\mathcal{F}))\la Ex(X,\mathcal{F})\la H^0(T^1(X,\mathcal{F})) \la H^2(\underline{Hom}(\Omega_X,\mathcal{F}))
\end{equation}
In our case, it follows from the above sequence that the obstruction to lift a deformation $H_B$ of $H$ over $B$ to a deformation $H_A$ over $A$, is in
$H^2(\underline{Hom}(\Omega_{{H_A}/A},\sheaf_H))=H^2(\underline{Hom}(\Omega_{{H_A}/A}\otimes \sheaf_H,\sheaf_H))=H^2(\underline{Hom}(\Omega_H,\sheaf_H))=0$. Hence $\phi$ is smooth.

Now assume that every irreducible component $C_i$ of $\Delta$ is a smooth rational curve and let $d_i = \deg( L_{{C_i}})$. If there is an $i$ such that $d_i <0$, then every section $s \in H^0(T^1_{qG}(H))$ 
vanishes along $C_i$ and hence every deformation of $H$ remains singular, and in fact not normal. Suppose that $d_i \geq 0$, $\forall i$. Since every component of $\Delta$ is rational, 
it is possible to find a section $s\in H^0(T^1_{qG}(H))$ such that it vanishes at $d_i$ points with order 1, and in fact we may choose these points to be normal crossing points. 
We will now show that locally the deformations of $H$ are unobstructed. This is clear for all index 1 points since they are complete intersections. 
It remains to consider slt points. So let $P\in H$ be an slt point. Then in suitable analytic coordinates this is isomorphic to  $(xy=0)/\mathbb{Z}_m(a,-a,1)$. 
Let $B \la A$ be a surjection of local Artin $\mathbb{C}$-algebras. Let $X_A$ be a 
$\mathbb{Q}$-Gorenstein deformation of $H$ over $A$. Let $p\colon \tilde{H} \la H$ be the index 1 cover of $H$. Then $X_A=\tilde{X_A}/\mathbb{Z}_m$, where $\tilde{X_A}$ is 
a deformation of $\tilde{H}$~\cite{KoBa88} over $A$. Now since $\tilde{H}$ is complete intersection, $\tilde{X_A}$ lifts to a deformation $\tilde{X_B}$ of $\tilde{H}$ over $B$ and hence 
$X_B=\tilde{X_B}/\mathbb{Z}_m$ is a deformation of $H$ over $B$. Therefore locally all $\mathbb{Q}$-Gorenstein deformations of $H$ are unobstructed. 
The obstruction to lift local deformations to global is in $H^2(T_H)$ which is zero by our assumptions. Hence if $d_i \geq 0$, $H$ is smoothable to a 3-fold $X$ with isolated cyclic quotient 
singularities and ordinary double points.

\end{proof}
Putting together Theorem~\ref{main-theorem} and Proposition~\ref{def-map} we get the following.
\begin{corollary}\label{global-def}
Let $\Delta \subset H$ be the germ of a surface along a proper curve $\Delta$ as in Theorem~\ref{main-theorem}.
In addition assume that $H^1(\sheaf_{\Delta})=0$, $H^2(T_H)=0$ and    \[
{\tilde{C}}^2 +p+2c_1+2c_2+\sum_{P\in U_3} \alpha_3(P) +\sum_{Q\in U_4} \alpha_4(Q) \geq 0\]
for every irreducible component $C$ of $\Delta$. Then there exists a one-parameter $\mathbb{Q}$-Gorenstein smoothing $X$ of $H$ whose singularities are
isolated cyclic quotient singularities of type $1/n(a,-a,1)$, or ordinary double points $(xy-zt=0)\subset \mathbb{C}^4$.
\end{corollary}

The following corollary shows that under certain conditions, it is possible to construct three dimensional terminal extremal neighborhoods $Y\la X$ such that $H\in |\sheaf_Y|$
\begin{corollary}\label{extremal}
Let $g \colon \Delta \subset H \la S \ni 0$ be a modification of an isolated surface singularity $0\in S$. Assume that $H$ has singularities of class $qG$ along $\Delta$
and that $R^1f_{\ast}\sheaf_H=0$.
Then a one parameter smoothing of $g$, $f\colon X \la Y$, with $X$ terminal, exists, if and only if \[
{\tilde{C}}^2 +p+2c_1+2c_2+\sum_{P\in U_3} \alpha_3(P) +\sum_{Q\in U_4} \alpha_4(Q) \geq 0\]
for every irreducible component $C$ of $\Delta$.
In particular, a terminal extremal neighborhood $f\colon \Delta \subset X \stackrel{f}{\la} Y \ni P$ exists such that $H \in |\sheaf_X|$, if and only if the above inequality holds
and $K_H \cdot C <0$, for every irreducible component of $\Delta$.
\end{corollary}
\begin{proof}
Since $R^1f_{\ast}\sheaf_H =0$ and the fibers of $f$ are at most 1-dimensional, it follows that $\Delta$ is a rational cycle of curves and $H^2(T_H)=0$. Hence by
Corollary~\ref{global-def} there exists a terminal smoothing $X$ of $H$. Then by~\cite{Wa76}, $g$ extends to $X$, and hence there exists a 3-fold contraction $g \colon X\la Y$, as claimed.
If it happens that $K_H \cdot C <0$, then $K_X \cdot C =K_H \cdot C <0$, and hence $Y$ is terminal as well and the contraction is an extremal neighborhood.

\end{proof}

\begin{example}
In this example we will construct a stable surface $Z$ such that locally around each singularity it has $\mathbb{Q}$-Gorenstein smoothings but 
it does not have any global $\mathbb{Q}$-Gorenstein smoothings. Such a surface will be in the boundary of the moduli space 
of stable surfaces.

Let $\sum_{i=1}^7E_i \subset U$ be the germ of a smooth surface around a chain of rational curves whose dual graph is 
\[
\xymatrix @R=-7pt @C=5pt{
                        &\underset{}{\overset{2}{\circ}} \ar@{-}[dr] &                                   &                    & \underset{}{\overset{2}{\circ}} \ar@{-}[dl]\\
                 &   &\underset{}{\overset{6}{\bullet}} \ar@{-}[r]\ar@{-}[dl]          & \underset{}{\overset{1}{\circ}} \ar@{-}[dr]   & \\
 \underset{}{\overset{2}{\circ}} \ar@{-}[r]  &   \underset{}{\overset{2}{\circ}}                                 &                                   &                    & \underset{}{\overset{3}{\circ}}
}
\]
where minus the selfintersection of a curve is denoted over it. Contracting all curves except the $-6$ and $-1$, we get a birational morphism $\Delta \subset \tilde{H} \stackrel{f}{\la} S\ni 0$, 
where $\Delta=C_1+C_2$, $C_i$ are the birational transforms of the $-6$ and $-1$ curves, and $0\in S$  is a quotient surface singularity whose dual graph is
\[
\overset{-2}{\circ}\mbox{\noindent ---}\overset{-2}{\circ}\mbox{\noindent ---}\overset{-3}{\circ}\mbox{\noindent ---}\overset{-2}{\circ}
\] 
and hence it is of type $1/11(1,8)$. The singular locus of $\tilde{H}$ consists of exactly four points $A_1,\; B_1 \in C_1$ and $A_2,\; B_2 \in C_2$ such that 
$A_1 \in \tilde{H} \cong 1/3(1,1)$, $A_2 \in \tilde{H} \cong 1/3(1,2)$, and both $B_1$, $B_2$ are of type $1/2(1,1)$. Now $C_1+C_2$ has an involution fixing the intersection $C_1 \cap C_2$ and 
interchanging $A_1$, $A_2$ and $B_1$, $B_2$. Identifying now $C_1$ and $C_2$ with the involution we get a surface germ $C\subset H$, $C\cong \mathbb{P}^1$, whose normalization is $\tilde{H}$. 
$H$ has semi-log-canonical singularities and in fact it has two semi-log-terminal singularities of types $(xy=0)/\mathbb{Z}_2(1,1,1)$, $(xy=0)/\mathbb{Z}_3(1,-1,1)$, one degenerate cusp 
of embedding dimension 3 and multiplicity 2, and is normal crossing elsewhere along $C$. Note that all the singularities of $H$ have $\mathbb{Q}$-Gorenstein smoothings. 
Let $\colon \tilde{H}\la H$ be the identification map. Straightforward calculations show that 
$K_H \cdot C = 1/6 > 0$ and that $\tilde{C}^2=(C_1+C_2)^2=-3$. Hence from Theorem~\ref{main-theorem} it follows that $\deg T^1_{qG}(H)\otimes \sheaf_C =-3+2=-1$, and hence 
$T^1_{qG}(H)\otimes \sheaf_C=\sheaf_{\mathbb{P}^1}(-1)$. Hence $\mathbb{T}^1_{qG}(H)=H^1(T_H)$, and therefore every $\mathbb{Q}$-Gorenstein deformation of $H$ is locally trivial. In particular, 
$H$ is not globally smoothable, even though it is locally smoothable. Notice also that in this case $H^2(T_H)=0$ which in the isolated singularities case is the obstruction space to lift 
local deformations to global ones. 

Next I claim that $H$ can be compactified to a stable surface $Z$. To do so we will first construct a terminal stable surface $W$ whith a singularity of type 
$1/11 (1,8)$ and then we will blow it up to get $Z$. So let $W \subset \mathbb{P}(1,1,11,8)$ be a general weighted hypersurface of degree 56, for example $x_0x_2^5+x_3^7+x_0^{56}+x_1^{56}=0$. 
This has exactly one singular point $P$ at $[0,0,1,0]$ which is of type $1/11(1,8)$. Moreover $K_W=\sheaf_W(56-1-1-11-8)=\sheaf_W(35)$ and hence ample. 

Now $U\la S$ is just a sequence of blow ups and hence $U\la S$ can be extended to a morphism $X \la W$. Contract all curves again except the $-6$ and $-1$ to get a compactification $\tilde{Z}$ 
of $\tilde{H}$, and identifying $C_1$, $C_2$ we get the desired compactification $Z$ of $H$. I now claim that $Z$ is stable. Indeed, this is true if and only if $\omega_Z^{[n]}$ is ample, for some $n$. 
This is true if and only if $\pi^{\ast}\omega_Z^{[n]}$ is ample, and by subadjunction, this is ample if and only $K_{\tilde{Z}}+C_1+C_2$ is ample. Moreover, let $g\colon \tilde{Z}\la W$ be the contraction map. 
Then $K_{\tilde{Z}}=g^{\ast}K_W+a_1C_1+a_2C_2$, with $a_1,\; a_2 \; >-1$, since the singularity of $W$ is log terminal. Hence 
\[
K_{\tilde{Z}}+C_1+C_2=g^{\ast}K_W+(a_1+1)C_1+(a_2+1)C_2
\]
Therefore $a_i+1 >0$, $i=1,2$, and hence if $\Gamma \subset \tilde{Z}$ is any curve different than $C_1$, $C_2$, then 
\[ K_{\tilde{Z}}\cdot \Gamma = 
K_W\cdot f(\Gamma) +(a_1+1)C_1\cdot \Gamma +(a_2+1)C_2\cdot \Gamma >0
\]
Moreover, a straightforward calculation shows that $ K_{\tilde{Z}}\cdot C_i =1/6>0$ and it now follows that $Z$ is stable.
\end{example}

\begin{example}
Let $U$ be the germ of a smooth surface around a chain of smooth
rational curves with the following dual graph\[
\stackrel{-2}{\circ} \mbox{\noindent ---} \stackrel{-2}{\circ} \mbox{\noindent ---}
\stackrel{-2}{\circ} \mbox{\noindent ---} \stackrel{-3}{\circ} \mbox{\noindent ---}
\underset{C_2}{\stackrel{-2}{\bullet}} \mbox{\noindent ---} \stackrel{-3}{\circ} \mbox{\noindent ---}
\underset{C_1}{\stackrel{-1}{\bullet}} \mbox{\noindent ---}
\stackrel{-2}{\circ} \mbox{\noindent ---} \stackrel{-5}{\circ}
\]
Contracting all curves except $C_1$ and $C_2$ we get a map $\tilde{H} \stackrel{f}{\la} S$, such that $0 \in S$ is
an $A_5$ singularity. Moreover, $\tilde{H}$ has exactly three singular points $P_1 \in C_1$, $P_2\in C_2$ and $Q \in C_1 \cap C_2$, analytically isomorphic to
$1/9 (1,5)$, $1/9(1,-5)$ and $1/3(1,1)$, respectively.
We can now identify $C_1$ and $C_2$ to a smooth rational curve $C$ and get a nonnormal surface $C\subset H$. $H$ has an \textit{slt} singular point analytically isomorphic to
$(xy=0)/\mathbb{Z}_9(5,-5,1)$, a degenerate cusp of type $T^3_{1,1}$, and is normal crossing everywhere else. Moreover, $\tilde{H}$ is the normalization
of $H$ and there is a morphism $g \colon H \la S$. An easy calculation now on $\tilde{H}$ shows that $(C_1+C_2)^2=-2/3$. Now from Theorem~\ref{main-theorem} it follows that
$d=\deg L= (C_1+C_2)^2+\alpha_3(Q)=-2/3+(1+1)/(1+1+1)=0$. Hence $L\cong \sheaf_C$, and hence there exists a one parameter terminal smoothing $Y$ of $H$ extending the local ones.
More precisely, $Y$ has exactly 1 singular point and it is of type $1/7(5,-5,1)$. The morphism $g$ extends to $Y$~\cite{Wa76}, and hence we get an extremal neighborhood
$f \colon C\subset Y \la X \ni 0$, such that $H\in |\sheaf_Y|$.

\end{example}
The next lemma describes the local deformations of singularities of class \textit{qG}.

\begin{lemma}
Let $(P\in H)$ be the germ of a singularity as in Theorem~\ref{main-theorem}. Then $(P\in H)$ has a smoothing $(P\in X)$ such that
\begin{enumerate}
\item $(P\in X)$ is smooth if $(P\in H)$ is not a degenerate cusp of embedding dimension 4 or an \textit{slt} point.
\item \[
(P\in X)\cong (xy-zt=0)\subset \mathbb{C}^4\]
if $(P\in H)$ is a degenerate cusp of embedding dimension 4.
\item \[
(P\in X) \cong \mathbb{C}^3/\mathbb{Z}_n(a,-a,1) \]
if $(P\in H)$ is an \textit{slt} point.
\end{enumerate}
\end{lemma}
\begin{proof}
The first part is obvious. Now suppose that $(P\in H)$ is a degenerate cusp of embedding dimension 4. Then in suitable local analytic coordinates it is given by equations
$xy-z^p-t^q=0$, $zt-x^r=0$. Then take $X$ to be $xy+zt-z^p-t^q-x^r=0$ which is analytically equivalent to $xy-zt=0$.

Let $(P\in H)$ be an \textit{slt} point. Then $(P\in H)\cong (xy=0)/\mathbb{Z}_n(a,-a,1)$. A $\mathbb{Q}$-Gorenstein smoothing is
$(xy+t=0)/\mathbb{Z}_n(a,-a,1,0) \cong \mathbb{C}^3/\mathbb{Z}_n(a,-a,1)$.
\end{proof}
The next lemma describes the support and the torsion part of $T^1_{qG}(H)$.
\begin{lemma}\label{torsion}
Let $\Delta \subset H$ be the germ of a surface along a proper curve as in Theorem~\ref{main-theorem}. Then
the support of $T^1_{qG}(H)$ is 1-dimensional. Its divisorial part is $\Delta$ and it has an embedded point over any pinch point or degenerate cusp
of type other than  $T^3_{\infty, \infty}$ and $T^4_{\infty,\infty,\infty}$. Moreover, let $C$ be any component of $\Delta$. Then \[
T^1_{qG}(H) \otimes \sheaf_C = L_C \oplus (\oplus_{i}\mathbb{C}_{{P_i}}) \]
where $L_C$ is a line bundle on $C$ and $\mathbb{C}_{{P_i}}$ is a 1-dimensional torsion sheaf supported at a degenerate cusp $P_i \in H$ of embedding dimension 4 that lies on $C$.
\end{lemma}
\begin{proof}
Locally, at an index 1 point, $H=(f=0)\subset \mathbb{C}^3$ and  \[
T^1(H)=T^1_{qG}(H)=\frac{\mathbb{C}[x,y,z]}{(f, J(f))},\]
where $J(f)$ is the Jacobian of $f$. Then at any generic point of $\Delta$, $H$ is normal crossing, i.e., $H=(xy=0)$ and therefore
$T^1_{qG}(H) = \mathbb{C}[x,y,z]/(x,y)=\mathbb{C}[z]$ and hence free of rank 1.

Let $P\in H$ be a pinch point. Then \[ T^1_{qG}(P\in H) =\frac{
\mathbb{C}[x,y,z]}{(x,y^2,yz)}\] and hence its support is defined by
the ideal $(x,y^2,yz)=(x,y^2,z)\cap (x,y)$. Therefore its divisorial
part is $\Delta$ and it also has an embedded point over $P$. Moreover,
$\Delta$ is irreducible at a neighborhood of $P$.

Let $P\in H$ be a degenerate cusp of type $T^2_{n}$, with $n \leq
\infty$. Then
\[ T^1_{qG}(P\in H) =\frac{
\mathbb{C}[x,y,z]}{(x,yz^2,(n+3)z^{n+2}-2zy^2)} \]
and hence its
support is defined by the ideal $I=(x,yz^2,(n+3)z^{n+2}-2zy^2)$ (if
$n =\infty $ we set as usual $z^n=0$). Therefore its divisorial part
is $\Delta$ and since $r(I:z^2)=(x,y,z)$, it also has an embedded
point over $P$(~\cite{AM69}).

Let $P\in H$ be a degenerate cusp of type $T^3_{p,q}$, with $p,\;q,\leq \infty$. Then
\[
 T^1_{qG}(P\in H) =\frac{\mathbb{C}[x,y,z]}{(px^{p-1}-yz,qy^{q-1}-xz,xy)}
\]
with the usual conventions about the cases that at least one of $p$, $q$, is infinity. Hence its support is defined by the ideal
$I=(px^{p-1}-yz,qy^{q-1}-xz,xy)$ and therefore its divisorial part is $\Delta$. If at least one of $p$, $q$ is not infinity, say $p$,
then $r(I:x)=(x,y,z)$ and hence the support of $T^1_{qG}(H)$ has an embedded point over $P$. If $p=q=\infty$, then
$I=(xy,xz,yz)$ and hence in this case there are no embedded points.

Let $C$ be any component of $\Delta$. Then in all of the above cases it is clear that
$T^1_{qG}(P\in H)\otimes \sheaf_{\Delta}$ is free of rank $1$.

Now let $P\in H$ be an \textit{slt} point. Let $P\in X$ be a $qG$-smoothing and let $\pi \colon \hat{X}\la X$ be the index 1 cover of $X$. Then $\pi^{-1}(H)=\hat{H} \la H$ is the index 1
cover of $H$~\cite{KoBa88} and hence $\hat{H}=(xy=0)\subset \mathbb{C}^3$. Therefore $qG$-smoothings of $H$ are quotients of $qG$-smoothings of $\hat{H}$ and therefore,
$T^1_{qG}(P\in H) =T^1(\hat{H})^{{\mathbb{Z}_n}}$ and hence $T^1_{qG}(P\in H)$ is torsion free and has support $\Delta$.

Finally, let $P\in H$ be a degenerate cusp of type $T^4_{p,q,r}$, with $p,\; q,\; r \leq \infty$. By Lemma~\ref{equations},
in suitable analytic coordinates, $P\in H$ is given by $(f(x,y,z,t)=g(x,y,z,t)=0)\subset \mathbb{C}^4$. We will now describe $T^1_{qG}(P\in H)$.
It is known that $T^1_{qG}(P\in H)=T^1(P\in H) =Ext^1_H(\Omega_H,\sheaf_H)$~\cite{Li-Sch67}. Moreover, there is an exact sequence \[
0 \la T_H \la T_{{\mathbb{C}^4}} \otimes \sheaf_H \la N_{{H/\mathbb{C}^4}} \la Ext^1_H(\Omega_H,\sheaf_H) \la 0\]
Writing down explicitely the isomorphisms $_{{\mathbb{C}^4}} \otimes \sheaf_H \cong \sheaf_H^{\oplus 4}$ and
$N_{{H/\mathbb{C}^4}}\cong \sheaf_H^{\oplus 2}$, it follows that $T^1_{qG}(P\in H)$ fits in the exact sequence \[
\sheaf_H^{\oplus 4} \stackrel{\phi}{\la} \sheaf_H^{\oplus 2} \la T^1_{qG}(P\in H)\la 0\]
where $\phi (1,0,0,0)=(\partial f/ \partial x, \partial g/ \partial x)$, $\phi (0,1,0,0)=(\partial f/ \partial y, \partial g/ \partial y)$,
$\phi (0,0,1,0)=(\partial f/ \partial z, \partial g/ \partial z)$ and $\phi (0,0,0,1)=(\partial f/ \partial t, \partial g/ \partial t)$.
Therefore we find that
\[
T^1_{qG}(P\in H) = \frac{\mathbb{C}[x,y,z,t] \oplus \mathbb{C}[x,y,z,t]}{M}
\]
where $M$ is the $\mathbb{C}[x,y,z,t]$-submodule of $\mathbb{C}[x,y,z,t] \oplus \mathbb{C}[x,y,z,t]$ that is generated by
$(f,0)$, $(0,f)$, $(g,0)$, $(0,g)$, $(\partial f/ \partial x, \partial g/ \partial x)$,
$(\partial f/ \partial y, \partial g/ \partial y)$, $(\partial f/ \partial z, \partial g/ \partial z)$ and
$(\partial f/ \partial t, \partial g/ \partial t)$. Now a straightforward calculation shows that the divisorial part of the support
of $T^1_{qG}(P\in H)$ is $\Delta$ and that if at least one of $p$, $q$, $r$ is not infinity, then it has an embedded point over $P$.
If $p=q=r=\infty$, then the support of $T^1_{qG}(P\in H)$ is defined by the ideal $(xy,xz,xt,yz,yt,zt)$ and it has no embedded point.
Moreover, a straightforward calculation shows that if $C$ is any irreducible component of $\Delta$, then
$T^1_{qG}(P\in H)\otimes \sheaf_{C} = \mathbb{C}_P \oplus \sheaf_{C}$
where $\mathbb{C}_P$ is a one-dimensional torsion sheaf concentrated at $P$ as claimed.
 \end{proof}

The next lemma shows how to calculate $T^1_{qG}(H)$ from a suitable embedding.

\begin{lemma}\label{degree-formula}
Let $\Delta \subset H$ be the germ of a surface along a proper curve as in the statement of Theorem~\ref{main-theorem}. 
Let $H \subset X$ be an embedding such that $X$ is a $\mathbb{Q}$-Gorenstein 3-fold,
 $H$ is Cartier in $X$, and for all $P\in \Delta \subset H$, $(P\in X)$ is a general $\mathbb{Q}$-Gorenstein smoothing of $(P\in H)$. 
Let $C$ be an irreducible component of $\Delta$ and $L_C$ be the free part of $T^1_{qG}(H)\otimes \sheaf_C$. Let $\phi$ be the natural map \[
\mathcal{N}_{H/X} \la \ext^1_H(\Omega_H,\sheaf_H) \]
Then  \[
L_C =\mathrm{Im} [ \mathcal{N}_{H/X}\otimes \sheaf_C \la \ext^1_H(\Omega_H,\sheaf_H)\otimes \sheaf_C] \]
and therefore \[
\deg (L_C \otimes \sheaf_C)= C \cdot H \]
\end{lemma}
\begin{proof}
Dualizing the exact sequence 
\[
0 \la I_{H,X}/I_{H,X}^2 \la \Omega_X \otimes \sheaf_H \la \Omega_H \la 0 \]
we get 
\begin{equation}
\cdots \la \mathcal{N}_{H/X} \stackrel{\phi}{\la}  \ext^1_H(\Omega_H,\sheaf_H)=T^1(H) \la  \ext^1_H(\Omega_X \otimes \sheaf_H ,\sheaf_H)\la 0
\end{equation}
The claim of the lemma then  can be checked locally. $\ext^1_H(\Omega_X \otimes \sheaf_H ,\sheaf_H)$ is torsion and it is supported over the \textit{slt} points of $H$ and the degenerate cusps of 
embedding dimension $4$. $X$ is smooth away from such points and Lemma~\ref{degree-formula}.1 is true around all points except possibly \textit{slt} points and degenerate cusps of embeding dimension 4. 
We will now check what happens around such points.

Let $P \in X$ be a singular point such that $P\in H$ is an \textit{slt} point. Let $P\in U$ be a small neighborhood of $X$ around $P$ and $V=U\cap H$. 
Then $(1)$ restricts to \[
 \sheaf_V \stackrel{\phi}{\la}  Ext^1_V(\Omega_V,\sheaf_V)
\]
Let $\pi \colon \hat{U} \la U$ be the index 1 cover of $U$. Then it is well known that $\hat{V}=\pi^{-1}(V)  \la V$ is the index 1 cover of $V$. Therefore this implies that \[
T^1_{qG}(V) = T^1(\hat{V})^{{\mathbb{Z}_n}}\]
Then $(1)$ follows from the diagram 
\[
\xymatrix{ 
\sheaf_{{\hat{V}}}^{{\mathbb{Z}_n}}  \ar@{=}[d]   \ar[r]^{{\hat{\phi}}} & T^1(\hat{V})^{{\mathbb{Z}_n}} \ar[d]\\
\sheaf_V \ar[r]^-{\phi} &   T^1(V)= Ext^1_V(\Omega_V,\sheaf_V) }
\] 
and from the fact that $\hat{\phi}$ is surjective. Hence there is a surjection \[
N_{V/U} \la T^1_{qG}(V) \]
It now remains to check what happens around a singular point $P\in H$ that is a degenerate cusp of embedding dimension 4. 
According to Lemma~\ref{torsion}, $T^1_{qG}(H)\otimes \sheaf_C = L_C \oplus \mathbb{C}_P$, where $\mathbb{C}_P$ is a one-dimensional torsion sheaf concentrated at $P$. 
The statement of the lemma will follow if we show that in a neighborhood of $P$, $ \mathrm{Ext}^1_H(\Omega_X \otimes \sheaf_H,\sheaf_H) \cong \mathbb{C}_P$.

Locally around $P$, $X=(xy-zw=0)\subset \mathbb{C}^4$. Hence there is an exact sequence 
\[
0 \la I_{X,{\mathbb{C}^5}}/I_{X,{\mathbb{C}^4}}^2 \la \Omega_{\mathbb{C}^5} \otimes \sheaf_X \la  \Omega_X \la 0 \]
Considering that $ I_{X,{\mathbb{C}^5}}/I_{X,{\mathbb{C}^5}}^2 \cong \sheaf_X$, this gives that 
\[
\mathrm{Hom}_H (\Omega_{{\mathbb{C}^5}}\otimes_{\sheaf_X} \sheaf_H,\sheaf_H) \la \mathrm{Hom}_H(\sheaf_H,\sheaf_H)\la 
\mathrm{Ext}^1_H(\Omega_X \otimes \sheaf_H,\sheaf_H) \la 0
\]
Now writing explicitely all the maps involved in the above sequence and the isomorphisms 
$\mathrm{Hom}_H (\Omega_{{\mathbb{C}^5}}\otimes_{\sheaf_X} \sheaf_H,\sheaf_H) \cong \sheaf_H^{\oplus 4}$, $\mathrm{Hom}_H(\sheaf_H,\sheaf_H)\cong \sheaf_H$, 
it follows that the last exact sequence is equivalent to 
\[
\sheaf_H^{\oplus 4} \stackrel{\phi}{\la} \sheaf_H \la  \mathrm{Ext}^1_H(\Omega_X \otimes \sheaf_H,\sheaf_H) \la 0
\]
where the map $\phi$ is given by sending the standard basis elements $(1,0,0,0)$, $(0,1,0,0)$, $(0,0,1,0)$ and $(0,0,0,1)$ to $x$, $y$, $z$, $t$, respectively. 
Hence it is now clear that $ \mathrm{Ext}^1_H(\Omega_X \otimes \sheaf_H,\sheaf_H) \cong \mathbb{C}_P$. This concludes the proof of the lemma.

\end{proof}

The next lemma was communicated to me by J\'anos Koll\'ar and it shows that it is always possible to find a 3-fold $X$ with the properties of Lemma~\ref{degree-formula}.

\begin{lemma}\label{patching}
Let $\Delta \subset H$ be the germ of a surface along a proper curve $\Delta$.
Assume that $H$ is smooth away from $\Delta$ and that generically along every component of $\Delta$ it has at worst singularities of type
$(xy=0)\subset \mathbb{C}^3$. Assume that for each $c \in \Delta$ we specify a neighborhood $U_c \subset X$ and a closed embedding $U_c \subset V_c$.
Assume that $V_c$ is smooth away from $c$ and singular for finitely many only $c$. Then there is a global embedding $\Delta \subset H \subset X$ that patches the local embedings.
\end{lemma}
\begin{proof}
Let $W=H-\{\mbox{a few transverse cuts}\}$. Then $W$ is Stein and it embedds in a $\mathbb{C}^n$, for some $n$. We want to get an embedding $W \subset V_W$ into something three dimensional
that has at most finitely many singularities along $\Delta$. Let $f \colon Y \la \mathbb{C}^n$ be the blow up of $\mathbb{C}^n$ along $W$. Let $E$ be the $f-$exceptional set. Then by the assumptions
on the singularities of $H$, a local calculation shows that away from finitely many singularities of $H$ that correspond to non normal crossing points, $E=E_1+E_2$, where $E_1$ and $E_2$ are
both smooth intersecting transversally. Moreover, $E_1+E_2$ is Cartier and the singular locus of $Y$ has codimension $3$ in $Y$. More precisely, locally the singularities of $Y$ are
of the type $(xy-zt=0)\subset \mathbb{C}^{n+1}$. Moreover, $\dim f^{-1}(c)=n-3$, for all $c\in \Delta$. Then by Bertini's theorem, the general $V^{\prime}=Z_1 \cap Z_2\cap \cdots \cap Z_{n-3}$,
with $Z_i \in |-E_1-E_2|$, has finitely many singularities and $V^{\prime} \la V_W=f_{\ast}(V^{\prime})$ is birational. $V_W$ has the required properties.

Therefore, by removing a few more transverse cuts, $ W=H-\{\mbox{a few transverse cuts}\}$ has a smooth embedding into a three dimensional space $V_W$.
Let $(H-W)\cap \Delta =\{c_1,\ldots, c_k\}$. Let $U_{{c_i}}$ be a Stein neighborhood of $c_i \in H$ and let $ U_{{c_i}}\subset V_{{c_i}}$ be an embedding such that $V_{{c_i}}$ is smooth
away from $c_i$. We want to glue $V_W$ and $V_{{c_i}}$ into a 3-fold $V$. This will follow from the following.

\noindent \textit{Claim:} Let $c\in \Delta \subset U$ be the germ of a surface singularity along a curve $\Delta$.
Assume that $U$ is singular along $\Delta$. Let $R\subset U$ be a small ring around $c\in U$. Let $R \subset X^1_R$, $R\subset X^2_R$ be two embeddings,
where $X^1_R$ and $X^2_R$ are smooth, which
we may also assume that they are Stein. Then, after shrinking $R$ if necessary, $X^1_R\cong X^2_R$.

The map $R \subset X^2_R$ extends to a map $X^1_R \la X^2_R$. Since $R$ is singular along $C\Delta$, then for all $c\in \Delta$, $T_c(R)=T_c(X^1_R)=T_c(X^2_R)$,
and hence by shrinking $R$ if necessary, it follows from the inverse function theorem that $X^1_R \cong X^2_R$.

Now gluing $V_W$ and $V_{{c_i}}$, $i=1,\ldots, k$, we obtain the required embedding.

\end{proof}
An immediate consequence of Lemmas~\ref{patching} and 3.6 is the following.
\begin{corollary}\label{embedding}
Let $\Delta \subset H$ be a surface germ as in Theorem~\ref{main-theorem}. Then there exists an embedding $H \subset X$ such that $X$ is a terminal $3$-fold and $H$ is Cartier in $X$.
Moreover, let $P\in X$ be a singular point
of $X$. Then either $P\in H$ is a degenerate cusp of embedding dimension 4 and $(P\in X) \cong (xy-zt=0)\subset \mathbb{C}^4$
or $P\in H$ is an \textit{slt} point, and  $(P \in X) \cong 1/n(a,-a,1)$.
\end{corollary}
We now want to calculate the degree of $L$, the free part of $T^1_{qG}(H)$ from the normalization of $H$. The next proposition is the key to do so.

\begin{proposition}\label{contraction}
Let $\Delta \subset X$ be the germ of a 3-fold along a curve $\Delta$ (proper or not). Suppose that $X$ has isolated singularities along $\Delta$ and that
locally in a neighborhood of any $P\in \Delta$, the general member of $|-K_X|$ that contains $\Delta$ has DuVal singularities.
Then there exists a divisorial contraction $E\subset Y \stackrel{f}{\la} X \supset \Delta$, such that $Y$ has canonical singularities and
\begin{enumerate}
\item $Y-E \cong X-\Delta$
\item $-K_Y$ is $f$-ample
\item Distinct irreducible components of $E$ are contracted onto distinct irreducible components of $\Delta$
\item The singularities of $Y$ lie on finitely many fibers over $X$
\end{enumerate}
Moreover, over the smooth locus of $X$, $f$ is the blow up of $X$ along $\Delta$.

In general it is hard to describe the singularities of $E$ and $Y$, but for the cases of interest in this paper we have the following.

\begin{enumerate}
\item Let $P$ be a point in $X$ such that $(P\in X)\cong (xy-zt=0)\subset \mathbb{C}^4$ and $\Delta$ is smooth at $P$.
Then $f^{-1}(P)$ is the union of two smooth rational curves $F_1$, $F_2$ intersecting
transversally at a point $Q$. $Y$ is smooth along $f^{-1}(P)-{Q}$ and  $(Q\in Y)\cong \mathbb{C}^3/\mathbb{Z}_2(1,1,1)$.
$E$ is smooth away from $Q$ and $(Q\in E)\cong (xy-z^2=0)\subset \mathbb{C}^3$. Moreover $K_Y \cdot F_i=-1/2$, $i=1,\;2$.
\item Let $P$ a point in $X$ such that $\Delta$ is smooth at $P$ and $(P\in \Delta \subset X)\cong (x=y=0) \subset \mathbb{C}^3/\mathbb{Z}_n(a,-a,1)$.
Then $(f^{-1}(P))_{red}=F=\mathbb{P}^1$. $E^{sing}\cap F=\{Q_1,Q_2\}$ such
that $(Q_1\in E)\cong \mathbb{C}^2/\mathbb{Z}_n(1,2a)$ and  $(Q_2\in E)\cong \mathbb{C}^2/\mathbb{Z}_n(1,-2a)$.
Moreover, $(Q_1 \in Y)\cong \mathbb{C}^3/\mathbb{Z}_n(a,-2a,1)$, $(Q_2 \in Y)\cong \mathbb{C}^3/\mathbb{Z}_n(-a,2a,1)$, $K_Y \cdot F=-1/n$, if $n$ is odd, and $K_Y\cdot F=-2/n$, if $n$ is even.
\end{enumerate}
Note that in the above two cases, there is always a DuVal $S\in |-K_X|$ that contains $\Delta$.
\end{proposition}
\begin{proof}

If a divisorial contraction $f\colon E \subset Y \la X \supset \Delta$ that satisfies the condititions of the proposition exists, then it is unique and in fact \[
Y\cong Proj \bigoplus_d i_{\ast} (I^d_{\Delta\cap U,U}) \]
where $U$ is the smooth part of $X$.
Therefore its existence is equivalent to the finite generation of the sheaf of algebras $ \oplus_d i_{\ast} (I^d_{\Delta\cap U,U})$~\cite{Tzi03}.
This can be checked locally and hence the existence of local
contractions imply the existence of a global one. We only need to check locally around the singular points of $X$.
Let $P\in X$ be a neighborhood of a singular point and let $S\in |-K_X|$ a DuVal surface that contains $\Delta$. Then by inversion of adjunction, $(X,S)$ is \textit{klt} and in fact
canonical since $K_X+S$ is Cartier. Let $g \colon Z \la X$ be the blow up of $X$ along $\Delta$. Let $E$ be the $g$-exceptional divisor that is over $\Delta$. Then
generically along  $\Delta$,
$K_Z=g^{\ast}K_X+E$ and $g^{\ast}S =S^{\prime}+E$. Hence $K_Z+S^{\prime}=g^{\ast}(K_X+S)$ and therefore $E$ is crepant for the pair $(X,S)$.
Hence there exists an extraction of $E$,
$f \colon E \subset \colon Y \la X \supset \Delta$, such that $f\colon Y-E \la X-\Delta$ is an isomorphism~\cite[Theorem 17.10]{Ko93}. 
Moreover, $K_Y+f_{\ast}^{-1}S=f^{\ast}(K_X+S)$ and hence $(Y,f_{\ast}^{-1}S)$
is canonical and therefore $Y$ has canonical singularities too.

Let $P\in X$ be a singular point such that $(P\in X)\cong (x=y=0)\subset \mathbb{C}^3/\mathbb{Z}_n(a,-a,1)$, where $x$, $y$, $z$ are the coordinates for $\mathbb{C}^3$. Hence $\Delta=\pi (L)$,
where $L$ is the line $x=y=0$. Then the divisorial contraction can be constructed from the following diagram.
\begin{equation}
\xymatrix{
B_L \hat{X}= \hat{Y}  \ar[d]_{{\hat{f}}}   \ar[r]^{\hat{\pi}}  & Y=\hat{X}/\mathbb{Z}_n \ar[d]^{f}\\
\mathbb{C}^3=\hat{X} \ar[r]^{\pi} &   X }
\end{equation}
where $\hat{Y}$ is the blow up of $\mathbb{C}^3$ along the $z$-axis, $x=y=0$. A straightforward calculation now of the lifting of the action of $\mathbb{Z}_n$ on $\hat{Y}$ shows that
the singular locus of $E$ is two points $Q_1$ and $Q_2$, and $(Q_1\in E)\cong 1/n(1,2a)$, $(Q_2\in E) \cong 1/n(1,-2a)$. Moreover, $(Q_1\in Y)\cong 1/n(a,-2a,1)$ and
$(Q_2\in Y)\cong 1/n(-a,2a,1)$. If $n$ is odd then the singular locus of $Y$ consists of the points $Q_1$ and $Q_2$. If on the other hand $n$ is even, then $Y$ may be singular along $f^{-1}(P)$.

It remains to find $K_Y\cdot F$, where $F=f^{-1}(P)_{red}$. Let $\hat{F}=\hat{f}^{-1}(0)$. Then $F=\hat{F}/\mathbb{Z}_n$. Let $u,\; u$ be homogeneous coordinates for $\mathbb{P}^1=\hat{F}$.
Let $\zeta$ be a primitive $n$-th root of unity. Then $\mathbb{Z}_n$ acts on $\hat{F}$ by $\zeta \cdot [u,v]=[\zeta^a u, \zeta^{-a}v]$. Hence $\hat{F}\la F$ is $n$-to-one if $n$ is odd and
$n/2$-to-one of $n$ is even. Now $\hat{\pi}$ is \'{e}tale in codimension 1 and hence $K_{{\hat{Y}}}=\pi^{\ast}K_Y$. Therefore \[
-1=K_{{\hat{Y}}}\cdot \hat{F}=\pi^{\ast}K_Y \cdot \hat{F}=K_Y \cdot \pi_{\ast}(\hat{F}) \]
and therefore $K_Y \cdot F =-1/n$, if $n$ is odd, and $K_Y \cdot F =-2/n$, if $n$ is even, as claimed.

Now let $P$ be a singular point of $X$ such that $(P\in X)\cong (xy-zt=0)\subset \mathbb{C}^4$. Then the divisorial contraction can be constructed by the following diagram~\cite{Tzi03}.
\begin{equation}
\xymatrix{
                              & Z \ar[dl]_{h} \ar[dr]^{\phi} &               \\
W \ar[dr]_{g}                  &                  &  Y \ar[dl]^{f}            \\
                           &     X              &     }
\end{equation}
where $W$ is the blow up of $X$ along $\Delta$. A straightforward calculation of the blow up shows that there are two $g$-exceptional divisors. A surface $E$ ruled over $\Delta$
and a $F\cong \mathbb{P}^2$
over $P$. The singular locus of $W$ consists of two points $R_1$ and $R_2$ on $E\cap F$. Then $Z$ is the blow up of $E$, $F_Z=h_{\ast}^{-1}F\cong F =\mathbb{P}^2$, and there are two
disjoint $h$-exceptional curves. Moreover, $\phi$ is the contraction
of $F_Z$. Calculations as in~\cite{Tzi03} show that $N_{{F_Z}/Z}=\sheaf_{{\mathbb{P}^2}}(-2)$ and hence it contracts to a $1/2(1,1,1)$ singularity. All other statements about the singularities
of $E$ can also be checked by local calculations.
\end{proof}

%% file: section3-tziolas.tex
\section{Proof of Theorem~\ref{main-theorem}}
\setcounter{equation}{0}

With assumptions and notation as in Theorem~\ref{main-theorem}, let $C$ be an irreducible component of $\Delta$, and $H \subset X$ an embedding as in Corollary 3.10.  Let 
$f \colon E \subset Y \la X \supset C$ be the divisorial contraction as in Proposition~\ref{contraction}, $H^{\prime}=f_{\ast}^{-1}H$ and $D^{\prime}=E\cdot H^{\prime}$. 
Let $\pi \colon \tilde{H} \la H$ be the normalization of $H$ and $\nu \colon \hat{H} \la H^{\prime}$ the normalization of $H^{\prime}$. 
Then there is a commutative diagram 
\begin{equation}
\xymatrix{ 
                              & \hat{H} \ar[dl]_{g} \ar[dr]^{\nu} &               \\
\tilde{H} \ar[dr]_{\pi}                  &                  &  H^{\prime} \ar[dl]^{f}            \\
                           &     H              &     }
\end{equation}

The next lemma describes the birational transform $H^{\prime}$ of $H$ in $Y$ 
and the minimal resolution of the normalization $\tilde{H}$ of $H$. 

\begin{lemma}\label{description-of-H}
\begin{enumerate}
\item Let $P\in H$ be a \textit{nc} or a pinch point. Then in a neighborhood of $P \in H$, $H^{\prime}$ is smooth and it is the normalization $\tilde{H}$ of $H$. Moreover,
$D^{\prime}=\tilde{C}$ is smooth, $P$ is not a ramification point of $\pi$, if $P$ is \textit{nc}, and it is a ramification point if it is a pinch point. 
\item Let $P\in H$ be a degenerate cusp with $\Gamma^2=-1$. Then $H^{\prime}=\tilde{H}$ is smooth. $D^{\prime}=\tilde{C}$ is the union of two smooth curves intersecting transversally, 
and hence globally $\tilde{C}$ is either two transversal curves or it is nodal. $P$ is a ramification of $\pi$.
\item Let $P\in H$ be a degenerate cusp with $\Gamma^2=-2$. Suppose it is of type $T^2_n$ with $n < \infty$. Then $H^{\prime}=\tilde{H}$. 
If it is of type $T^2_{\infty}$, then $H^{\prime}$ is normal crossing and its normalization is $\tilde{H}$, which is smooth. 
In both cases $D^{\prime}$ is the union of two smooth curves intersecting transversally 
at a point $Q$, $(Q \in H^{\prime}) \cong (xy-z^n=0)$ in the $T^2_n$ case and nc in the the $T^2_{\infty}$. Moreover,in both cases $\tilde{C}$ is 
either two transversal curves or nodal and  $P$ is a ramification of $\pi$.
\item Let $P\in H$ be a degenerate cusp with $\Gamma^2=-3$. Then $D^{\prime}=E\cdot H^{\prime} =C^{\prime}+F$, 
where $F\cong \pone $ is $f$-exceptional and $C^{\prime}$ is smooth and $2-1$ over $C$. Moreover, $C^{\prime}\cap F=\{P,\; Q\}$, with $P\neq Q$ and,
\begin{enumerate}
\item Suppose $P \in H$ is of type $T^3_{p,q}$, for $1 \leq p,\; q < \infty$. Then $H^{\prime}$ is normal and is the blow up of $H$ along $C$. 
$P\in H^{\prime}$ is an $A_{p-1}$ DuVal singularity and $Q\in H^{\prime}$ is an $A_{q-1}$ DuVal singularity. Moreover, $F^2=-1-1/p-1/q$, $C^{\prime}$ 
is of type $FC_1$ in both $(P\in H^{\prime})$ and $(Q\in H^{\prime})$.
\item Suppose $P \in H$ is of type $T^3_{p,\infty}$. Then $P \in H^{\prime}$ is an $A_{p-1}$ singularity and $Q \in H^{\prime}$ is nc. 
Let $\hat{F}$ be the divisorial part of $\nu^{-1}(F)$. Then $\hat{F}^2=-1-1/p$.
\item Suppose $P \in H$ is of type $T^3_{\infty,\infty}$. Then $H^{\prime}$ is nc and hence its normalization $\hat{H}$ is smooth. Moreover, $P$, $Q\in H$ are nc points and 
$\hat{F}^2=-1$.
\end{enumerate} 
\item Let $P\in H$ be a degenerate cusp with $\Gamma^2=-4$. Then it is of type $T^4_{p,q,r}$, with $2 \leq p,\; q,\; r \leq \infty$. Let $U$ be the minimal resolution of $\tilde{H}$.
\begin{enumerate}
\item Suppose $ p,\; q,\; r < \infty$. Then $H^{\prime}$ is normal. Then the extended dual graph of $C^{\prime} \subset \tilde{H}$ is
\begin{enumerate}
\item If $r \geq 3$, $D$, it is 
\[
\underset{{C^{\prime}}}{\bullet}\mbox{\noindent ---}\overset{-2}{\underset{{E_{p-2}}}{\circ}}\mbox{\noindent ---} \cdots \mbox{\noindent ---}
\overset{-2}{\underset{E_1}{\circ}}\mbox{\noindent ---}
\overset{-3}{\underset{{F_1}}{\circ}} 
\mbox{\noindent ---}\overset{-2}{\underset{{D_1}}{\circ}}\mbox{\noindent ---}
\cdots \mbox{\noindent ---}\overset{-2}{\underset{{D_{r-3}}}{\circ}} \mbox{\noindent ---}\overset{-3}{\underset{{F_2}}{\circ}} 
\mbox{\noindent ---}\overset{-2}{\underset{B_1}{\circ}}\mbox{\noindent ---}\cdots \mbox{\noindent ---}
\overset{-2}{\underset{{{B_{q-2}}}}{\circ}}\mbox{\noindent ---}\underset{{C^{\prime}}}{\bullet}
\]
\item If $r=2$, it is \[
\underset{{C^{\prime}}}{\bullet}\mbox{\noindent ---}\overset{-2}{\underset{{E_{p-2}}}{\circ}}\mbox{\noindent ---} \cdots \mbox{\noindent ---}
\overset{-2}{\underset{E_1}{\circ}}\mbox{\noindent ---}\overset{-4}{\underset{F}{\circ}} 
\mbox{\noindent ---}\overset{-2}{\underset{{B_1}}{\circ}}\mbox{\noindent ---}
\cdots \mbox{\noindent ---}\overset{-2}{\underset{{B_{q-2}}}{\circ}} \mbox{\noindent ---}\underset{{C^{\prime}}}{\bullet} 
\]
\end{enumerate}
In both cases $H^{\prime}$ is obtained from $U$ by contracting all curves except $E_1$ and $B_1$. If $p$ or $q$ is $2$ then we set $E_i=\emptyset$ or $B_i=\emptyset$. 
In particular, if $p=q=2$ then $H^{\prime}=\tilde{H}$. 
$P$ is a ramification point of $\hat{H} \la H$, and $\tilde{C}$ is either a union of two smooth curves intersecting transversally, or nodal. 
Moreover, $D^{\prime}=C^{\prime}+E_1^{\prime}+B_1^{\prime}$, where $E_1^{\prime}$ and $B_1^{\prime}$ are the birational transforms 
of $E_1$ and $B_1$ in $H^{\prime}$ and 
\begin{gather*}
{E_1^{\prime}}^2=-1-\frac{1}{p-2}+\frac{2r-3}{4r-4} \\
{B_1^{\prime}}^2=-1-\frac{1}{q-2}+\frac{2r-3}{4r-4}
\end{gather*}
\item Suppose $p,\; q < \infty$ and $r =\infty$. Then the singular locus $\Delta$ of $H$ is the transversal union of $C$ and another smooth curve $C_1$ and,
\begin{enumerate}
\item The extended dual graph of $\Delta \subset \tilde{H}$ is
\[
\underset{{C^{\prime}}}{\bullet}\mbox{\noindent ---}\overset{-2}{\underset{{E_{p-2}}}{\circ}}\mbox{\noindent ---} \cdots \mbox{\noindent ---}
\overset{-2}{\underset{E_1}{\circ}}\mbox{\noindent ---}\overset{-2}{\underset{F_1}{\circ}}\mbox{\noindent ---}\underset{{C_1^{\prime}}}{\bullet}
\mbox{\noindent ---}\overset{-2}{\underset{F_2}{\circ}}
\mbox{\noindent ---}\overset{-2}{\underset{{B_1}}{\circ}}\mbox{\noindent ---}
\cdots \mbox{\noindent ---}\overset{-2}{\underset{{B_{q-2}}}{\circ}} \mbox{\noindent ---}\underset{{C^{\prime}}}{\bullet} 
\]
\item $\hat{H}$ is obtained from $U$ by contracting all curves except $E_1$ and $B_1$ and 
$D^{\prime}=C^{\prime}+E_1^{\prime}+B_1^{\prime}$,  where $E_1^{\prime}$ and $B_1^{\prime}$ are the birational transforms 
of $E_1$ and $B_1$ in $H^{\prime}$. Moreover, let $\hat{E_1}$ and $\hat{{B_1}}$ be the birational transforms of $E_1$ and $B_1$ in $\hat{H}$. Then
\begin{gather*}
{\hat{E_1}}^2=-\frac{p}{2(p-2)} \\
{\hat{B_1}}^2=-\frac{q}{2(q-2)}
\end{gather*}
\end{enumerate}
\item Suppose $p,\; r  < \infty$ and $q=\infty$. Then the singular locus of $H$, $\Delta$, is the transversal union of two smooth curves $C_1$, and $C_2$. 
Then up to an analytic change of coordinates we can assume that $H=(xy-z^p=zt-x^r=0)\subset \mathbb{C}^4$ and $C=(x=z=t=0)$.
Then the extended dual graph of $\Delta \subset \tilde{H}$ is
\[
\underset{{C_1^{\prime}}}{\bullet}\mbox{\noindent ---}\overset{-2}{\underset{{B_{p-2}}}{\circ}}\mbox{\noindent ---} \cdots \mbox{\noindent ---}
\overset{-2}{\underset{B_1}{\circ}}\mbox{\noindent ---}\overset{-3}{\underset{F}{\circ}}
\mbox{\noindent ---}\overset{-2}{\underset{{E_{1}}}{\circ}}\mbox{\noindent ---} \cdots \mbox{\noindent ---}
\overset{-2}{\underset{E_{r-2}}{\circ}}
\mbox{\noindent ---}\underset{{C_2^{\prime}}}{\bullet} \mbox{\noindent ---}\underset{{C_1^{\prime}}}{\bullet}
\]
and the extended dual graph of $\hat{H}$ is
\[
\underset{{C_1^{\prime}}}{\bullet}\mbox{\noindent ---}\overset{-2}{\underset{{B_{p-2}}}{\circ}}\mbox{\noindent ---} \cdots \mbox{\noindent ---}
\overset{-2}{\underset{B_1}{\circ}}\mbox{\noindent ---}\overset{-3}{\underset{F}{\circ}}
\mbox{\noindent ---}\overset{-2}{\underset{{E_{1}}}{\circ}}\mbox{\noindent ---} \cdots \mbox{\noindent ---}
\overset{-2}{\underset{E_{r-3}}{\circ}}\mbox{\noindent ---}
\overset{-3}{\underset{E_{r-2}}{\circ}}\mbox{\noindent ---}
\overset{-1}{\underset{F_1}{\circ}}
\mbox{\noindent ---}\underset{{C_2^{\prime}}}{\bullet} \mbox{\noindent ---}\underset{{C_1^{\prime}}}{\bullet}
\]
Moreover, $\hat{H}$ is obtained from $U$ by contracting all curves except $B_1$ and $F_1$ and 
$D^{\prime}=C^{\prime}+E_1^{\prime}+F_1^{\prime}$,  where $E_1^{\prime}$ and $F_1^{\prime}$ are the birational transforms 
of $F_1$ and $B_1$ in $H^{\prime}$. Moreover, let $\hat{F_1}$ and $\hat{{B_1}}$ be the birational transforms of $F_1$ and $B_1$ in $\hat{H}$. Then
\begin{gather*}
{\hat{B_1}}^2=-1-\frac{1}{p-2}+\frac{2r-3}{4r-4} \\
{\hat{F_1}}^2=-1+\frac{2r-3}{4r-4}
\end{gather*}

\item Suppose $p < \infty$ and $q=r =\infty$. Then the singular locus of $H$, $\Delta$, is the transversal union of three smooth curves $C_1$, $C_2$ and $C_3$. Moreover,
the extended dual graph of $\Delta \subset \tilde{H}$ is
\[
\underset{{C_1^{\prime}}}{\bullet}\mbox{\noindent ---}\overset{-2}{\underset{{E_{p-2}}}{\circ}}\mbox{\noindent ---} \cdots \mbox{\noindent ---}
\overset{-2}{\underset{E_1}{\circ}}\mbox{\noindent ---}\overset{-2}{\underset{F}{\circ}}\mbox{\noindent ---}\underset{{C_2^{\prime}}}{\bullet}
\mbox{\noindent ---}\underset{{C_3^{\prime}}}{\bullet} \mbox{\noindent ---}\underset{{C_1^{\prime}}}{\bullet}
\]
and the extended dual graph of $\hat{H}$ is
\begin{enumerate}
\item If $C$ is $C_1$ or $C_2$, it is
\[
\underset{{C_1^{\prime}}}{\bullet}\mbox{\noindent ---}\overset{-2}{\underset{{E_{p-2}}}{\circ}}\mbox{\noindent ---} \cdots \mbox{\noindent ---}
\overset{-2}{\underset{E_1}{\circ}}\mbox{\noindent ---}\overset{-2}{\underset{F}{\circ}}\mbox{\noindent ---}\underset{{C_2^{\prime}}}{\bullet}
\mbox{\noindent ---}\overset{-2}{\underset{F_1}{\circ}}\mbox{\noindent ---}\overset{-1}{\underset{B}{\circ}}
\mbox{\noindent ---}\underset{{C_3^{\prime}}}{\bullet} \mbox{\noindent ---}\underset{{C_1^{\prime}}}{\bullet}
\]
In which case  $\hat{H}$ is obtained by contracting all curves except $E_1$ and $B$ and 
$D^{\prime}=C^{\prime}+E_1^{\prime}+B^{\prime}$,  where $E_1^{\prime}$ and $B^{\prime}$ are the birational transforms 
of $E_1$ and $B$ in $H^{\prime}$. Moreover, let $\hat{{E_1}}$ and $\hat{{B}}$ be the birational transforms of $E_1$ and $B$ in $\hat{H}$. Then
\begin{gather*}
{\hat{E_1}}^2=-\frac{p}{2(p-2)} \\
{\hat{B_1}}^2=-\frac{1}{2}
\end{gather*}
\item If $C=C_3$, then it is 
\[
\underset{{C_1^{\prime}}}{\bullet}\mbox{\noindent ---}\overset{-1}{\underset{B_1}{\circ}}
\mbox{\noindent ---}\overset{-3}{\underset{{E_{p-2}}}{\circ}}\mbox{\noindent ---}\overset{-2}{\underset{{E_{p-3}}}{\circ}}
\mbox{\noindent ---} \cdots \mbox{\noindent ---}
\overset{-2}{\underset{E_1}{\circ}}\mbox{\noindent ---}\overset{-3}{\underset{F}{\circ}}\mbox{\noindent ---}\overset{-1}{\underset{B_2}{\circ}}
\mbox{\noindent ---}\underset{{C_2^{\prime}}}{\bullet}
\mbox{\noindent ---}\underset{{C_3^{\prime}}}{\bullet} \mbox{\noindent ---}\underset{{C_1^{\prime}}}{\bullet}
\]
and in this case  $\hat{H}$ is obtained by contracting all curves except $B_1$ and $B_2$ and 
$D^{\prime}=C^{\prime}+B_1^{\prime}+B_2^{\prime}$,  where $B_1^{\prime}$ and $B_2^{\prime}$ are the birational transforms 
of $B_1$ and $B_2$ in $H^{\prime}$. Moreover, let $\hat{{B_1}}$ and $\hat{{B_2}}$ be the birational transforms of $B_1$ and $B_2$ in $\hat{H}$. Then
\[
{\hat{B_1}}^2={\hat{B_2}}^2=-1+\frac{2p-3}{4p-4}
\]
\end{enumerate}
\item Suppose that $p=q=r=\infty$. Then $\tilde{H}$ is smooth. Moreover, $D^{\prime}=E_1+E_2+C^{\prime}$, where the $E_i$ are $f$-exceptional, and if $\hat{E}_i$ 
are the divisorial parts of $\nu^{-1}(E_i)$, then $\hat{E}_i^2=-1/2$. 
\end{enumerate}
\item If $P\in H$ is an \textit{slt} point, then $H^{\prime}=\tilde{H}$, $D^{\prime}=\tilde{C}$ 
and $f^{-1}(P)=\{R_1,R_2\}$ where $ (R_1 \in \tilde{H})\cong 1/n(1,a)$ and $(R_2 \in \tilde{H})\cong 1/n(1,-a)$.
\end{enumerate}
\end{lemma}
\begin{remark}
In the case that any of $p$, $q$, $r$ is infinity, then the selfintersection numbers of the exceptional curves that appear in the previous lemma are 
exactly the corresponding limits to infinity of the selfintersection numbers that were found in the finite case.
\end{remark}

\begin{proof}
Let $P\in X$ be a point such that $P\in H$ is not a degenerate cusp of embedding dimension 4 or an \textit{slt} point. Then $X$ is smooth at $P$ and therefore, $Y$ is just the blow up 
of $X$ along $C$ and hence $H^{\prime}$ is also the blow up of $H$ along $C$. Then an explicit calculation of the blow up using the normal forms for the equations of the singularities that are given 
in Lemma~\ref{equations}, show parts (1)-(4).   

Let $P\in H$ be an \textit{slt} point. Then $(P\in H)\cong \hat{H} /\mathbb{Z}_n(a,-a,1)$, where $\hat{H}=(xy=0)\subset \mathbb{C}^3$. Let $\hat{U}$ be the blow up of $\hat{H}$ along 
$\hat{C}=(x=y=0)$. Then from the proof of Proposition~\ref{contraction} and in particular diagram (4), it follows that $\hat{H}=\pi^{-1}(H)$, $H^{\prime}=\hat{U}/\mathbb{Z}_n$ and hence 
it is the normalization of $H$, $\tilde{H}$. The claim about the singularities of $\tilde{H}$ follows from~\cite{KoBa88}, or by a careful calculation in diagram (4) in Proposition~\ref{contraction}.

It only remains to show part $(5)$ of the lemma. I will only indicate how to do the first part when $p,\; q,\; r<\infty $. The rest are similar. 
In this case, $P\in H$ is a degenerate cusp of embedding dimension 4 of type $T^4_{p,q,r}$ with  $p,\; q,\; r<\infty $ and the singular locus $\Delta$ of $H$ is $C$. First 
we describe the minimal resolution of the normalization of $H$. 
Let $f \colon U \la H$ be the minimal semiresolution, $\pi \colon U^{\prime} \la U$ the normalization of $U$, $\Gamma$ the reduced part of $f^{-1}(P)$ and $C$ the dual curve of $U$. 
Let $E_1, \ldots , E_n$ 
the $f$-exceptional curves, $C^{\prime}=\pi^{-1}(C)$ and $F_i=\pi^{\ast}(E_i)$. Then by subadjunction, $\omega_{{U^{\prime}}}=\pi^{\ast} \omega_U \otimes \sheaf_{{U^{\prime}}}(-C^{\prime})$, 
we find that \[
E_i \cdot \Gamma = F_i^2 +2 -F_i \cdot C^{\prime}\]
Now considering that $\Gamma$ is a cycle and $\Gamma \cdot C =2$, it follows that
\[
\Gamma^2 =\sum_{i=1}^n F_i^2 +2n-2 \]
In our case $\Gamma^2=-4$ and hence \[
\sum_{i=1}^nF_i^2=-2-2n\]
and hence, considering that $\tilde{H}$ has cyclic quotient singularities and the extended dual graph of $\tilde{H}$ is a chain~\cite{Ba83}, 
we get the two possibilities for the fundamental cycle that are claimed in $5.a$. Now one easily sees that \[
(C^{\prime}+F_1+F_2+\sum E_i +\sum D_i +\sum B_i)\cdot E =0 \]
where $E$ is any of the $E_i$, $B_i$ or $D_i$. Hence the blow up of $H$ along $C$ is obtained from $U^{\prime}$ be contracting all exceptional curves except those of self intersection $\leq -3$. 
An explicit calculation of the blow up now using the normal forms given in Lemma~\ref{equations}, shows the relation between the $p,\; q,\; r$ and the number of the $-2$ exceptional curves. 

Alternatively, one could first describe explicitely the blow up $g\colon H_W \la H$ of $H$ along $C$ by using the normal forms given in Lemma~\ref{equations}. 
Straightforward calculations show that $g^{-1}(P)=F_1+F_2$ and $g^{-1}(C) = C^{\prime}+F_1+F_2$. Moreover, $C^{\prime}\cap F_i =P_i$, $P_1 \notin F_2$, 
$P_2\notin F_1$. $P_1 \in H_W$ is an $A_{p-2}$ DuVal double point and $P_2 \in H_W$ is an $A_{q-2}$ DuVal double point. Now one can easily find the minimal resolution 
of $\tilde{H}$. The abstract method described above has the advantage that it allows us directly to find all the possible 
fundamental cycles for the singularities in question without explicit calculations using normal forms. That would be particularly useful in dealing with cusps 
of higher embedding dimension. 

We now want to describe $H^{\prime}$. According to diagram $(5)$ in the proof of Proposition~\ref{contraction}, $H^{\prime}$ fits in the commutative diagram
\begin{equation*}
\xymatrix{ 
                              & H_Z \ar[dl]_{h} \ar[dr]^{\phi} &               \\
H_W \ar[dr]_{g}                  &                  &  H^{\prime} \ar[dl]^{f}            \\
                           &     H              &     }
\end{equation*}
In the notation and setting of diagram $(5)$, $H_W=g_{\ast}^{-1}H$ is the blow up of $H$ along $C$.  Hence $E \cap H_W = C^{\prime}$,  and $F\cap H_W =F_1+F_2$,
as we saw earlier. Moreover, $F_i \not \subset E$, $i=1,\; 2$. Now $Z$ is the blow up of $W$ along $E$ and hence 
$H_Z=h_{\ast}^{-1}H_W$ is the blow up of $H_W$ along $C^{\prime}$. A similar argument as above shows that the $h$-exception curves are $E_1$ and $B_1$. 
Finally $\phi$ contracts $F=\mathbb{P}^2$, and 
hence $H^{\prime}$ is obtained from $H_Z$ by contracting $F_1$ and $F_2$ and the $f$-exceptional curves correspond to $E_1$ and $B_1$ as claimed. Now, explicit calculations 
using the description of the extended dual graphs, give the claims about ${E_1^{\prime}}^2$ and ${B_1^{\prime}}^2$.
\end{proof}

Let $L$ be the locally free part of $T^1_{qG}\otimes \sheaf_C=T^1_{qG}(H)/(torsion)$ and $d=\deg L$.
Let $H \subset X$ be an embedding as in Corollary~\ref{embedding}, and let $f\colon E \subset Y \la X \supset C$ the divisorial contraction as in Proposition~\ref{contraction}. 
Then $D^{\prime}=H^{\prime} \cdot E =C^{\prime}+F_0$, where $F_0$ is supported on finitely many fibers which by 
Lemma~\ref{description-of-H} must be over degenerate cusps with $\Gamma^2=-3$ and $\Gamma^2=-4$ of type $T^4_{p,q,r}$ with $p,\; q,\; r \leq \infty$ and 
at least one of $p,\; q$ is $\geq 3$. 
More precisely, let $P\in H$ be a degenerate cusp of type $T^4_{p,q,r}$. Then $f^{-1}(P)=F_1+F_2$, $F_i \cong \mathbb{P}^1$ and $Q\in E$ is an $A_1$ double point, where $Q =F_1\cap F_2$. 
Then in a neighborhood of $P\in H$, $F_0=f^{-1}(P)=F_1+F_2$ and $Q \not \in C^{\prime}$, if $p,\; q \geq 3$,  and $F_0=F_1$ if one of $p,\; q$ is 2. 
In the last case, $C^{\prime}$ goes through the singular point of $E$. 
Over a point $P\in H$ of type $T^4_{2,2,r}$ or such that $\Gamma^2=-1,\; -2$, $H^{\prime}=\tilde{H}$ and $D^{\prime}=C^{\prime}$ either breaks into the union of two 
transversal curves or it has a node, as shown by Lemma~\ref{description-of-H}.  
$E$ is normal, its singularities are over \textit{slt} points of $H$ and over degenerate cusps of 
embedding dimension 4 and they are described explicitly in Proposition~\ref{contraction}. Moreover from the description of $H^{\prime}$ given in Lemma~\ref{description-of-H} it is clear 
that $C^{\prime}$ goes through exactly those singularities of $E$ that are over \textit{slt} points of $H$ or over degenerate cusps of type $T^4_{2,q,r}$. 

It must also 
be pointed out that since $E$ is $\mathbb{Q}$-Cartier in $Y$, $E\cdot H^{\prime}=D^{\prime}$ is also $\mathbb{Q}$-Cartier in $H^{\prime}$ and hence intersection numbers make sense there. 
However, in the degenerate cases when at least one of $p, \; q, \; r$ is $\infty$, $H^{\prime}$ is not normal and $C^{\prime}$ is not $\mathbb{Q}$-Cartier in $H^{\prime}$.

Since $H$ is singular along $C$ and generically \textit{nc}, $f$ is generically the blow up of $C$ and hence 
\begin{equation}
f^{\ast}H=H^{\prime}+2E
\end{equation}
Moreover, $E^2=-\delta + F_1$, where $\delta$ is a section over $C$ and $F_1$ is supported on finitely many fibers. Intersecting (1) with $E$ we find 
\begin{equation}
f^{\ast}(D)=H^{\prime}\cdot E +2E^2=D^{\prime}-2\delta+2F_1
\end{equation}
where $D\in Pic(C)$ and $\deg(D)=H \cdot C =d$. Now it is possible to write $K_E=aD^{\prime}+f^{\ast}(D_1)$, for some number $a$ and $D_1 \in Pic(C) \otimes \mathbb{Q}$. Intersecting 
with a general fiber $l$ we find that $a=-1$ and hence 
\begin{equation}
K_E=-D^{\prime}+f^{\ast}(D_1)
\end{equation}
Intersecting now with $D^{\prime}$ we find that 
\begin{equation}
K_E \cdot D^{\prime} +{D^{\prime}}^2 = 2\deg (D_1)
\end{equation}
where all intersection numbers above are in $E$. Taking into account that $D^{\prime}=C^{\prime}+F_0$ we find that \[
K_E \cdot D^{\prime} +{D^{\prime}}^2=K_E \cdot C^{\prime} +{C^{\prime}}^2 +(K_E\cdot F_0 + 2C^{\prime}\cdot F_0)-\frac{\mu}{2}\]
where $\mu$ is the number of singularities $P\in H$ of type $T^4_{2,q,r}$, with $q \geq 3$. Over a neighborhood of a singularity of this type, $f^{-1}(P)=F_1+F_2$, and $F_0=F_1$. 
Moreover, $F_i^2=-1/2$, and thus the term $-\mu /2$ in the above formula.

\textit{Claim:} \[
K_E\cdot F_0 + 2C^{\prime}\cdot F_0=2c_3+2c_4-2m\]
where $c_3,\; c_4$ are the  number of degenerate cusps with $\Gamma^2=-3,\; -4$ and $m$ the number of degenerate cusps of type $T^4_{2,2,r}$ that lie on $C$.

Let $P\in H$ be a degenerate cusp with $\Gamma^2=-3, -4$. Then from the earlier discussion, $f^{-1}(P)=F_0=\mathbb{P}^1$ if $\Gamma^2=-3$ and $f^{-1}(P)=F_1+F_2$, if $\Gamma^2=-4$. 
Moreover, $F_0=F_1+F_2$ if $P\in H$ is of type $T^4_{p,q,r}$, $p,q\geq 3$, and $F_0=F_1$, if $p=2$, $q\geq 3$. 
Suppose that $P \in H$ is of type $T^4_{p,q,r}$, $p,q\geq 3$. Then $K_E\cdot F_0+2C^{\prime}\cdot F_0=2$. 

Now suppose that $P\in H$ is of type $T^4_{2,q,r}$, $q \geq 3$. Then \[
K_E\cdot F_0 +2F_0 \cdot C^{\prime}=K_E\cdot F_1+2F_1\cdot C^{\prime}=-1+2(1/2+1)=2\]
as before. A similar calculation for points with $\Gamma^2=-3$ shows the claim.

Therefore (5) becomes 
\begin{equation}
K_E \cdot C^{\prime} +{C^{\prime}}^2 +2c_3+2c_4-2m-\frac{\mu}{2}= 2\deg (D_1)
\end{equation}
By adjunction 
\[
K_E \cdot C^{\prime} +{C^{\prime}}^2=2p_a(C^{\prime})-2+\mathrm{Diff}(C^{\prime},E) 
\]
where $\mathrm{Diff}(C^{\prime},E)$ is the different of $C^{\prime}$ in $E$. Now (5) becomes 
\begin{equation}
\deg(D_1)=p_a(C^{\prime})-1+c_3+c_4-m-\frac{\mu}{4}+\frac{1}{2} \mathrm{Diff}(C^{\prime},E)
\end{equation}
Intersecting (3) with $\delta$ in $E$ we find that \[
K_E \cdot \delta = -D^{\prime} \cdot \delta +f^{\ast}(D_1)\cdot \delta = -D^{\prime}\cdot \delta +\deg (D_1)\]
Intersecting (2) with $\delta$ in $E$ we find \[
d=D^{\prime}\cdot \delta -2\delta^2+2F_1\cdot \delta\]
where $E^2=-\delta +F_1$. Taking into account this, (6) becomes 
\begin{equation}
K_E\cdot \delta +\delta^2 =-d-\delta^2+2F_1\cdot \delta +p_a(C^{\prime})-1+c_3+c_4-m-\frac{\mu}{4}+\frac{1}{2}\mathrm{Diff}(C^{\prime},E) 
\end{equation}
\textit{Claim:} 
\begin{equation}
-\delta^2 +2\delta \cdot F_1 =\frac{1}{2}{D^{\prime}_{{H^{\prime}}}}^2+\frac{1}{2}d
\end{equation}
where by ${D^{\prime}_{{H^{\prime}}}}^2$ we mean the self intersection of $D^{\prime}$ in $H^{\prime}$.

Intersecting (2) with $E$ and taking into account that $E\cdot \delta = -\delta^2+\delta \cdot F_1$ and ${D^{\prime}_{{H^{\prime}}}}^2=E\cdot D^{\prime}$,  we find that 
\[
-d=E\cdot D^{\prime} -2(E\cdot \delta -E\cdot F_1) = {D^{\prime}_{{H^{\prime}}}}^2-2(E\cdot \delta -E \cdot F_1) \]
But $E\cdot F_1=(-\delta+F_1)\cdot F_1 =-\delta \cdot F_1$. Hence $E\cdot \delta -E\cdot F_1 =-\delta^2+2\delta \cdot F_1$. Therefore \[
-d={D^{\prime}_{{H^{\prime}}}}^2-2(-\delta^2+2\delta \cdot F_1)\]
and hence \[
-\delta^2 +2\delta \cdot F_1 =\frac{1}{2}{D^{\prime}_{{H^{\prime}}}}^2+\frac{1}{2}d\]
as claimed. Taking this into account, (7) becomes 
\begin{equation}
\delta \cdot K_E +\delta^2 = -d +\frac{1}{2}{D^{\prime}_{{H^{\prime}}}}^2+\frac{1}{2}d+p_a(C^{\prime})-1+c_3+c_4-m-\frac{\mu}{4}+\frac{1}{2}\mathrm{Diff}(C^{\prime},E)
\end{equation}
Let $\mathrm{Diff}(\delta,E)$ be the different of $\delta$ in $E$. Then by adjunction \[
K_E\cdot \delta +\delta^2=2p_a(C)-2+\mathrm{Diff}(\delta, E) \]
and (9) becomes 
\begin{equation}
d={D^{\prime}_{{H^{\prime}}}}^2+2(p_a(C^{\prime})-1)-2(2p_a(C)-2)+2(c_3+c_4-m)-\frac{\mu}{2}+\mathrm{Diff}(C^{\prime},E)-2\mathrm{Diff}(\delta,E)
\end{equation}
We now want to calculate the differents $\mathrm{Diff}(C^{\prime},E)$ and $\mathrm{Diff}(\delta,E)$. The next lemma allows us to do so. 

\begin{lemma}
Let $0\in C \subset H$ be the germ of a cyclic quotient singularity of type $\frac{1}{n}(1,a)$ around a smooth curve $C$. Let $f\colon U \la H$ be the minimal 
resolution. Assume that the extended dual graph is \[
\underset{C}{\bullet}\mbox{\noindent ---}\underset{{E_{1}}}{\circ}\mbox{\noindent ---} \cdots \mbox{\noindent ---}\underset{E_m}{\circ}
\]
Then\[
\mathrm{Diff}(C,0\in H)=1-1/n\] 
\end{lemma}
\begin{proof}
In the minimal resolution $U$ we have the following equations
\begin{gather*}
K_U=f^{\ast}K_H+\frac{a-n+1}{n}E_1+\{ other\} \\
f^{\ast}C=C_U+\frac{a}{n}E_1+\{ other\}
\end{gather*}
Hence \[
2p_a(C)-2=K_U\cdot C_U +C_U^2=K_H\cdot C+C^2 -1+1/n \]
and the lemma follows.

\end{proof}

From Proposition~\ref{contraction} and Lemma~\ref{description-of-H} it follows that the singular points of $E$ lie over \textit{slt} points of $H$ and 
degenerate cusps with $\Gamma^2=-4$. More precisely, let $P\in H$ be an \textit{slt} point. Then $E$ has two singular points on $f^{-1}(P)$ which are of type $1/n(1,2a)$ and $1/n(1,-2a)$. 
Let $P \in H$ be a degenerate cusp with $\Gamma^2=-4$. Then $f^{-1}(P)=F_1+F_2$ and $E$ has an $A_1$ double point at $Q=F_1\cap F_2$. 

$C^{\prime}$ goes through exactly those singular points of $E$ that lie over \textit{slt} points of $H$ and degenerate cusps of type $T^4_{2,q,r}$. 
Over $T^4_{2,2,r}$ points, $C^{\prime}$ is either nodal or normal crossing at the singular point $Q\in E$ and it is easy to see that in such a case, $C^{\prime}$ is 
Cartier in $E$ and hence these points contribute nothing to the different $\mathrm{Diff}(C^{\prime},E)$. 

Let $s$ and $\mu$ be the numbers of \textit{slt} and $T^4_{2,q,r}$, $q\geq 3$, points respectively. Then by the previous lemma 
\begin{gather*}
\mathrm{Diff}(C^{\prime},E)=\sum_{i=1}^s\left( \mathrm{Diff}(C^{\prime},P_i \in E)+ \mathrm{Diff}(C^{\prime},P_{-i} \in E)\right)+
\sum_{Q} \mathrm{Diff}(C^{\prime},Q\in E)=\\
2(s-\sum_{i=1}^s\frac{1}{n_i})+\frac{\mu}{2}
\end{gather*}
where $(P_i \in E)\cong 1/{n_i}(1,2a_i)$, $(P_{-i}\in E)\cong 1/{n_i}(1,-2a_i)$ are the singular points over \textit{slt} points 
and $(Q\in E)\cong (xy-z^2=0)\subset \mathbb{C}^3$ are the singularities over points of type $T^4_{2,q,r}$, $q\geq 3$.

Now $\delta$ is a section and it goes through one of the singular points of $E$ that lie over \textit{slt} points of $H$, but it also goes through the singular point of $E$ that lies over 
a degenerate cusp of embedding dimension 4. The last assertion can either be seen by direct calculations from Proposition~\ref{contraction}, or as follows. Let $P\in H$ be a degenerate 
cusp of embedding dimension 4 and let $l=f^{-1}(0)=l_1+l_2$. By Proposition~\ref{contraction}, $E$ is singular at $Q= l_1\cap l_2$. Then if $\delta$ does not go through this singularity, 
then $-1/2=K_Y \cdot l_i =E \cdot l_i =(-\delta +l)\cdot l_i =\delta \cdot l_i \in \mathbb{Z}$. Therefore to calculate the different of $\delta$ in $E$ it is necessary to also 
consider the singular points that lie over degenerate cusps of embedding dimension 4. By Proposition~\ref{contraction}, all such points are isomorphic to $(xy-z^2=0)\subset \mathbb{C}^3$. 
Then by applying the previous lemma once more we find $\mathrm{Diff}(\delta,Q\in E)=1-1/2=1/2$, and hence  
\[
\mathrm{Diff}(\delta,E)=\sum_{i=1}^{{c_{4}}} \mathrm{Diff}(\delta,Q_i\in E) +\sum_{i=1}^s \mathrm{Diff}(\delta,R_i\in E)=
s-\sum_{i=1}^s\frac{1}{n_i} +\frac{c_{4}}{2}
\]
where $Q_i\in E$ are the singular points over the degenerate cusps of $H$ and $R_i\in E$ is $P_i$ or $P_{-i}$ as above over the \textit{slt} points. Therefore we find that 
\[
\mathrm{Diff}(C^{\prime},E)-2\mathrm{Diff}(\delta,E)=-c_4+\frac{\mu}{2} \]
Hence (10) becomes
\begin{equation}
d={D^{\prime}_{{H^{\prime}}}}^2+2(p_a(C^{\prime})-1)-2(2p_a(C)-2)+2c_3+c_4-2m
\end{equation}
The next simple lemma gives a singular version of the ramification formula for curves.
\begin{lemma}
Let $\pi \colon C^{\prime} \la C$ be a finite morphism between two curves with $C$ smooth. Suppose that $C^{\prime}$ has at worst nodes as singularities and let $n$ be the number of them. 
Then \[
2p_a(C^{\prime})-2= \deg \pi (2p_a(C)-2)+\deg R +2n \]
where $R$ is the part of the ramification divisor in the smooth part of $C^{\prime}$. 
\end{lemma}
To see the proof of it, let $p: \tilde{C}\la C^{\prime}$ be the normalization of $C^{\prime}$ (i.e, the blow up of the nodes). Then use the standard ramification formula for $\pi \circ p$ 
and the fact that $p_a(C^{\prime})=p_a(\tilde{C})+n$.

In our case now, $C^{\prime}$ is ramified exactly over the pinch points of $H$ and over degenerate cusps with $\Gamma^2=-1, -2$ as well as over points of type $T^4_{2,2,r}$. 
Moreover, $C^{\prime}$ is smooth over the pinch points but it is nodal or nc over the degenerate cusps. Hence by the previous lemma
\[
2p_a(C^{\prime})-2=2(2p_a(C)-2)+p+2c_1+2c_2+2m \]
where $p$, $c_1$ and $c_2$ are the numbers of pinch points and degenerate cusps with $\Gamma^2=-1,\; -2$ respectively. Then (11) becomes 
\begin{equation}
d={D^{\prime}_{{H^{\prime}}}}^2+p+2c_1+2c_2+2c_3+c_4= {\hat{D}}^2+p+2c_1+2c_2+2c_3+c_4
\end{equation}
where $\hat{D}=\nu^{\ast}(D^{\prime})$. Let $\hat{C}$ be the divisorial part of $\nu^{-1}(C^{\prime})$.
Next we will get a formula for $d$ that involves $\hat{C}^2$ instead of $\hat{D}^2$. Then we will reduce it to the normalization $\tilde{H}$ of $H$. 
Let $U_3$ be the set of points $P$ such that $P\in H$ is a degenerate cusp with $\Gamma^2=-3$, $W_4$ the set of points $Q$ such that $Q \in H$ is a degenerate cusp of type $T^4_{p,q,r}$ 
with $p,q\geq 3$, and $V_4$ the set of points $R$ such that $R \in H$ is a degenerate cusp of type $T^4_{2,q,r}$ with $q \geq 3$. 

Let $\alpha $ be a divisor in $H^{\prime}$. In what follows we will denote by $\hat{\alpha}$ the divisorial part of $\nu^{-1}(\alpha)$. 

By Proposition~\ref{contraction} 
and Lemma~\ref{description-of-H}, if $P\in U_3$, then near $P$, $F_0=f^{-1}(P)=F_P\cong \mathbb{P}^1$, 
if $Q\in W_4$, then near $Q$, $F_0=f^{-1}(Q)=F_{1,Q}+F_{2,Q}$, where $F_{i,Q}$ are smooth rational curves 
intersecting transversally at a point. Finally if $R\in V_4$, then $f^{-1}(Q)=F_{1,R}+F_{2,R}$ and near $R$, $F_0=F_R$, where $F_R=F_{i,R}$ for $i=1$ or $2$.
Hence,
\begin{equation}
\hat{D}=\hat{C}+\sum_{{P\in U_3}}\hat{F}_P +\sum_{{Q\in W_4}}(\hat{F}_{1,Q}+\hat{F}_{2,Q}) +\sum_{R \in V_4} \hat{F}_R
\end{equation}
Therefore
\begin{equation}
\hat{D}^2=\hat{C}^2+\sum_{{P\in U_3}} \beta_3(P) + \sum_{{Q\in W_4}}\beta_4(Q) +\sum_{R \in V_4}\beta^{\prime}(R) 
\end{equation}
where $\beta_3(P)=\hat{F}_P^2+2\hat{C}\cdot \hat{F}_P$, $\beta_4(Q)=(\hat{F}_{1,Q}+\hat{F}_{2,Q})^2 + 2 \hat{C}\cdot (\hat{F}_{1,Q}+\hat{F}_{2,Q})$ and 
$\beta_4^{\prime}(R)=\hat{F}_R^2+2\hat{C}\cdot \hat{F}_R$. 

Let $P\in U_3$, $Q\in W_4$ and $R\in V_4$. Then $P\in H$ is of type $T^3_{p,q}$, $Q\in H$ of type $T^4_{p,q,r}$ and $R\in H$ of type $T^4_{2,q,r}$. Then 
according to Lemma~\ref{description-of-H}.4 and some straightforward calculations in the minimal resolution of $\hat{H}$ we find that, 
\begin{gather}
\hat{F}_P^2=-1-\frac{1}{p}-\frac{1}{q}\\
\hat{F}_{1,Q}^2=-1-\frac{1}{p-2}+\frac{2r-3}{4r-4}\\
\hat{F}_{2,Q}^2=-1-\frac{1}{q-2}+\frac{2r-3}{4r-4}\\
\hat{F}_{1,Q}\cdot \hat{F}_{2,Q} =\frac{1}{4r-4}\\
\hat{C}\cdot \hat{F}_{1,Q}=\frac{1}{p-2}\\
\hat{C}\cdot \hat{F}_{2,Q}=\frac{1}{q-2}
\end{gather}
and therefore,
\begin{gather}
\beta_3(P)=-1+\frac{1}{p}+\frac{1}{q}\\
\beta_4(Q)=-1 +\frac{1}{p-2}+\frac{1}{q-2}\\
\beta_4^{\prime}(R)=\frac{1}{q-2}-\frac{2r-3}{4r-4}
\end{gather}
where in the case that any of the $p$, $q$, $r$ is infinity, we understand the above formulas to mean the corresponding limits at infinity. 

Let $g\colon \hat{H} \la \tilde{H}$ be the contraction of the $f$-exceptional curves, i.e., of $\hat{F}_0$. Then 
\begin{equation}
g^{\ast} \tilde{C} = \hat{C}+\sum_{P\in U_3} \gamma_P \hat{F}_P +\sum_{Q\in W_4} (\gamma_{1,Q}\hat{F}_{1,Q}+\gamma_{2,Q} \hat{F}_{2,Q}) +\sum_{R \in V_4}\gamma_R^{\prime} \hat{F}_R
\end{equation}
and therefore
\begin{equation}
{\tilde{C}}^2=\hat{C}^2+\sum_{P\in U_3} \delta_3(P) +\sum_{Q\in W_4} \delta_4(Q) +\sum_{R\in V_4}\delta^{\prime}_4(R)
\end{equation}
where $\delta_3(P)=(\gamma_P \hat{F}_P)\cdot \hat{C}$, $\delta_4(Q)=(\gamma_{1,Q}\hat{F}_{1,Q}+\gamma_{2,Q}\hat{F}_{2,Q})\cdot \hat{C}$ and 
$\delta^{\prime}_4(R)=\gamma^{\prime}_R\hat{F}_R \cdot \hat{C}$. 

\textit{Claim:} 
\begin{enumerate}
\item Let $P\in U_3$ be a degenerate cusp of type $T^3_{p,q}$. Then
\[
\delta_3(P)=\frac{(p+q)^2}{pq(p+q+pq)}\]
\item Let $Q \in W_4$ be a degenerate cusp of type $T^4_{p,q,r}$ with. Then 
\[
\delta_4(Q)=\frac{1}{p-2}+\frac{1}{q-2}-\frac{r(p+q)-4}{rpq-p-q}\]
\item Let $R\in V_4$ be a degenerate cusp of type $T^4_{2,q,r}$. Then \[
\delta^{\prime}_4(R)=\frac{(4r+q-6)^2}{4(2rq-q-2)(r-1)(q-2)}
\]
\end{enumerate}
where again in the case that any of $p$, $q$, $r$ is infinity, we understand the above formulas to mean the corresponding limits at infinity.

Assuming the claim, then from (13), (15) and (24) it follows that 
\begin{equation}
d={\tilde{C}}^2+\sum_{P\in U_3} (\beta_3(P)-\delta_3(P)) +\sum_{Q\in W_4}(\beta_4(Q)-\delta_4(Q)) +\sum_{R\in V_4}(\beta^{\prime}_4(R)-\delta^{\prime}_4(R)) +k
\end{equation}
where $k=p+2c_1+2c_2+2c_3+c_4$. On the other hand,
\begin{gather*}
\beta_3(P)-\delta_3(P)=-1+\frac{1}{p}+\frac{1}{q}-\frac{(p+q)^2}{pq(p+q+pq)}=-1+\frac{p+q}{p+q+pq}\\
\beta_4(Q)-\delta_4(Q)= -1+\frac{1}{p-2}+\frac{1}{q-2}- \frac{1}{p-2}-\frac{1}{q-2}+\frac{r(p+q)-4}{rpq-p-q}=-1+\frac{r(p+q)-4}{rpq-p-q}\\
\beta^{\prime}_4(R)-\delta^{\prime}_4(R)=\frac{1}{q-2}-\frac{2r-3}{4r-4} - \frac{(4r+q-6)^2}{4(2rq-q-2)(r-1)(q-2)}= -1+\frac{r(2+q)-4}{2rq-q-2}
\end{gather*}
Hence \[
d={\tilde{C}}^2 +p+2c_1+2c_2+\sum_{P\in U_3} \alpha_3(P) +\sum_{Q\in U_4} \alpha_4(Q) \]
where $U_3$ and $U_4$ are the sets of degenerate cusps of multiplicity 3 and 4, respectively. Let $P\in U_3$ be of type $T^3_{p,q}$. Then  \[
\alpha_3(P)=1+\frac{p+q}{pq+p+q}\] 
let $Q\in U_4$ be of type $T^4_{p,q,r}$. Then \[
\alpha_4(Q)=\frac{r(p+q)-4}{rpq-p-q}\]
This concludes the proof of Theorem~\ref{main-theorem}.

The claim follows from a lengthy but straightforward calculation. I will indicate how to do the second part about $\delta_4(Q)$. This is local over a degenerate cusp $Q\in H$ with $\Gamma^2=-4$. 
Suppose that the cusp is of type $T^4_{p,q,r}$ with $r\geq 3$. Then locally around $Q$ we can write \[
  g^{\ast}\tilde{C} =\hat{C}+\gamma_{1,Q}F_{1,Q}+\gamma_{2,Q}F_{2,Q}\]
Intersecting with $F_{1,Q}$ and $F_{2,Q}$ and taking into consideration (16)-(21), it follows that 
\begin{gather*}
4(r-1)-2(pr+p-4)\gamma_{1,Q}+(p-2)\gamma_{2,Q}=0\\
4(r-1)+(q-2)\gamma_{1,Q}-2(qr+q-4)\gamma_{2,Q}=0
\end{gather*}
A straightforward calculation now using (18)-(24) shows that \[
\delta_4(Q)=\frac{1}{p-2}\gamma_{1,Q}+\frac{1}{q-2}\gamma_{2,Q}=\frac{1}{p-2}+\frac{1}{q-2}-\frac{r(p+q)-4}{rpq-p-q}
\]